\documentclass[12pt]{amsart}

\usepackage[left=3cm,top=2.5cm,bottom=2.5cm,right=3cm]{geometry}
\usepackage{times}
\usepackage{amssymb}
\usepackage{graphicx,xspace}
\usepackage{epsfig}
\usepackage{enumitem}
\usepackage[usenames,dvipsnames]{xcolor}
\usepackage{tikz}
\usepackage{tikz-cd}
\usepackage{gensymb}
\usepackage[T1]{fontenc}
\usepackage[utf8]{inputenc}
\usepackage{bm}
\usepackage[boxsize=1em]{ytableau}
\usepackage[presets={abc,vec-cev}]{letterswitharrows}
\usepackage[config, labelfont={normalsize}]{caption,subfig}
\usepackage{mathrsfs}
\definecolor{green}{RGB}{0,127,0}
\definecolor{red}{RGB}{191,0,0}
\usepackage[colorlinks,cite color=red,link color=green,pagebackref=true]{hyperref}
\usepackage{todonotes}

\tikzset{
	arrow/.pic={\draw[every arrow/.try,->,>=#1] (0,0) -- +(.1pt,0);},
	pics/arrow/.default={triangle 90}
}

\newcommand{\Part}{\operatorname{Part}}

\usetikzlibrary{matrix,arrows,calc}
\usetikzlibrary{decorations.markings}
\usetikzlibrary{patterns}
\usetikzlibrary{snakes}
\usetikzlibrary{shapes}

\usepackage[capitalize]{cleveref}

\theoremstyle{plain}

\newtheorem{lemma}{Lemma}[section]
\newtheorem{theorem}[lemma]{Theorem}
\newtheorem{corollary}[lemma]{Corollary}
\newtheorem{proposition}[lemma]{Proposition}

\theoremstyle{remark}
\newtheorem{remark}[lemma]{Remark}

\newtheorem{definition}[lemma]{Definition}

\newcommand{\svdots}{%
	\vbox{
		\scriptsize \baselineskip 2pt \lineskiplimit 0pt
		\hbox {.}\hbox {.}\hbox {.}\kern-0.75pt
	}%
}

\newcommand{\shdots}{%
	\hbox{
		\scriptsize \baselineskip 20pt \lineskiplimit 2pt
		\hbox {.}\hbox {.}\hbox {.}\kern2.95pt
	}%
}

\DeclareMathOperator{\LLT}{LLT}

\newcommand{\Z}{\mathbb{Z}}

\def\la{\lambda}

\def\ka{\kappa}

\newcommand{\lla}{\bm{\la}}
\newcommand{\mmu}{\bm{\mu}}
\newcommand{\nnu}{\bm{\nu}}

\DeclareMathOperator{\SYT}{SYT}
\DeclareMathOperator{\SSYT}{SSYT}

\DeclareMathOperator{\inv}{inv}
\DeclareMathOperator{\Inv}{Inv}

\author[M.~Kowalski]{Maciej Kowalski}
\address{
Institute of Mathematics, 
Polish Academy of Sciences, 
ul. Śniadeckich 8, 
00-956 Warszawa, Poland}
\email{mkowalski@impan.pl}

\title [A combinatorial formula for LLT cumulants of melting lollipops]{A combinatorial formula for LLT cumulants of melting lollipops in terms of spanning trees}

\begin{document}

\maketitle

\begin{abstract}
	We prove a combinatorial formula for LLT cumulants of melting lollipops as a positive combination of LLT polynomials indexed by spanning trees. Due to Schur-positivity of LLT polynomials, the result implies an affirmative answer to a question from \cite{DolegaKowalski2021} for this class of unicellular LLT cumulants, and gives an independent proof of their Schur-positivity. In the special case of the complete graph, we also express the formula in terms of parking functions.
\end{abstract}

\ytableausetup{smalltableaux}

\section{Introduction}

\textit{LLT polynomials} were introduced by Lascoux, Leclerc, and Thibon~\cite{LascouxLeclercThibon1997} in the context of quantum groups. Their importance in the theory of symmetric functions quickly became apparent mostly due to the rich yet hidden combinatorics and a broad spectrum of relations --- for instance, with the Kazhdan--Lusztig theory \cite{GrojnowskiHaiman2007}, algebraic geometry \cite{LeclercThibon2000}, and knots theory \cite{Haglund2016}.

The motivation for this paper originated from another closely related field: the celebrated Macdonald polynomials. In 2005, Haglund, Haiman, and Loehr \cite{HaglundHaimanLoehr2005} noticed that Macdonald polynomials naturally decompose as a positive combination of LLT polynomials. This, in particular, meant that Macdonald's conjecture about Schur-positivity of Macdonald polynomials \cite{Macdonald1988} (proved by Haiman \cite{Haiman2001}) could get a new proof by showing a similar result for their LLT counterparts. This was indeed achieved by Grojnowski and Haiman \cite{GrojnowskiHaiman2007} using the aforementioned Kazhdan--Lusztig theory. In fact, their result covers the most general case of arbitrary skew shapes, thus extending the results from \cite{LeclercThibon2000} for the special case of tuples of straight shapes.

In 2017, Dołęga introduced \emph{Macdonald cumulants} \cite{DolegaFeray2017,Dolega2019}: roughly speaking, the higher order approximation of the product of Macdonald polynomials as $q\rightarrow 1$. In particular, the cumulants of order one are Macdonald polynomials themselves. Using the introduced terminology, Dołęga conjectured that Macdonald cumulants are Schur-positive \cite{Dolega2017}. The proof would generalize the celebrated result of Haiman \cite{Haiman2001} and might shed new light on the connection between symmetric functions, representation theory, and algebraic geometry. Broadlyly speaking, Haiman was able to relate the coefficients in the Schur expansion of Macdonald polynomials to Hilbert schemes and to the characters of the symmetric group's action on a certain Hecke algebra. Thus, by the definition of a character, Haiman proved that the Schur polynomial coefficients are positive integers. Therefore, since Macdonald cumulants generalize Macdonald polynomials, Dołęga's conjecture suggests the existence of a class of geometric objects or representations that the cumulants correspond to in a similar sense. Both of these observations could be potential approaches to the proof.

In \cite{DolegaKowalski2021}, the author together with Dołęga introduced \emph{LLT cumulants}: the LLT equivalent of Macdonald cumulants. What is more, it turns out that an analogous property to the one from \cite{HaglundHaimanLoehr2005} holds for cumulants: Macdonald cumulants expand as a positive combination of LLT cumulants \cite[Theorem 2.6]{DolegaKowalski2021}. Therefore, it proved natural to conjecture Schur-positivity of LLT cumulants \cite[Conjecture 2.5]{DolegaKowalski2021} as a means to achieve Schur-positivity of Macdonald cumulants.

In the special case of unicellular LLT cumulants, based on extensive computer simulations, we have noticed that it is always possible to find an expansion of the following special form:
\begin{equation} \label{conj:Tdecgen}
\ka(\lla/\mmu) = \sum\limits_T \LLT(\nu(T)),
\end{equation}
where $\lla/\mmu$ is a sequence of unicellular shapes and the sum runs over some spanning trees of $G_{\lla/\mmu}$ and $\nu(T)$ denotes a sequence of vertical-strip shapes corresponding to $T$ (see \cref{sec:lltgraphs} and \cref{sec:cumunicell} for details). Even though we were unable to devise a clear pattern in the general case, we managed to find one when $\lla/\mmu$ corresponds to a melting lollipop, and it is tempting to believe that such a pattern exists for all unicellular LLT cumulants.

\begin{definition} \label{def:meltinglollipop}
	For $k,l,m\in\mathbb{Z}_{\ge 0}$, the \emph{melting lollipop} $L^{(k)}_{(l,m)}$ is a graph on $l+m$ vertices obtained by joining the complete graph $K_m$ to a path of length $l$ at a vertex $v\in K_m$ and erasing $k$ edges of $K_m$ that are incident to $v$.
\end{definition}

In our setting, melting lollipops appear as \emph{LLT graphs}. Heuristically speaking, the vertices of such a melting lollipop correspond to the cells of a sequence of unicellular shapes and the edges correspond to pairs which can contribute an inversion in an LLT polynomial (see \cref{sec:preliminaries} and \cref{sec:meltinglollipops} for the precise definitions). In this case, we have managed to find a pattern that transforms \eqref{conj:Tdecgen} into a surprisingly simple combinatorial formula. This gives our main result.

\begin{theorem} \label{thm:meltinglollipopgeneralization}
	Let $\lla/\mmu$ be a sequence of $l+m$ unicellular diagrams corresponding to a melting lollipop $L^{(k)}_{(l,m)}$. Then
	\begin{equation} \label{eq:meltingtheorem}
	\ka(\lla/\mmu) = \sum_{T\subseteq L^{(k)}_{(l,m)}} \LLT(\nu(T)),
	\end{equation}
	where the sum runs over all spanning trees of $L^{(k)}_{(l,m)}$.
\end{theorem}

Note that thanks to the result of Grojnowski and Haiman \cite{GrojnowskiHaiman2007}, as an immediate consequence of the above theorem, we get Schur-positivity of LLT cumulants in the special case of $\lla/\mmu$ corresponding to a melting lollipop. It is an independent proof to the one presented in \cite{DolegaKowalski2021}.

Moreover, if we apply the result to the case $k=l=0$, i.e., when the melting lollipop is a complete graph, we can further express the formula in terms of parking functions.

\begin{corollary} \label{thm:singcelldecomp}
	Let $\lla/\mmu = ((1),\dots,(1))$ be a sequence of $m$ unicellular non-skew shapes. Then,
	\begin{equation} \label{thm:formulawithtrees}
	\ka(\lla/\mmu) = \sum_{T\subseteq K_m} \LLT(\nu(T)) = \sum_{f\in PF_{m-1}} \LLT(\mu(f)),
	\end{equation}
	where the first sum runs over all Cayley trees with $m$ vertices, the second over all parking functions on $m-1$ cars, and $\nu(T)$ and $\mu(f)$ are sequences of vertical-strip shapes associated to such objects (see \cref{sec:cumunicell} and \cref{sec:schroderpaths}).
\end{corollary}

We would like to mention the similarity between \eqref{thm:formulawithtrees} and the celebrated formula for the character of the diagonal coinvariant algebra (conjectured in \cite{HaglundHaimanLoehrRemmelUlyanov2005} and proved in \cite{CarlssonMellit2018}). Both formulas are given as linear combinations of LLT polynomials associated to parking functions, but the LLT polynomials corresponding to a fixed parking function differ in the two. It would be interesting to find a representation-theoretical interpretation of LLT cumulants, which encourages further investigation.

The basic approach for understanding LLT cumulants is to switch from the standard framework, where the combinatorial objects of focus are Young diagrams, Dyck, or Schr\"oder paths, to a more general one: that of graphs with colored vertices, which we call \emph{LLT graphs} (originally defined in \cite{DolegaKowalski2021}). The author together with Dołęga applied this reasoning to study the $e$-positivity phenomenon that gained significant attention in the recent years and, as a result, were able to expand on the existing theory (see, e.g., \cite{Alexandersson2021, AlexanderssonSulzgruber2020}).

In this paper, the focus on graphs corresponding to sequences of partitions allows us to transform the generating functions corresponding to the LHS of \eqref{thm:formulawithtrees} and \eqref{eq:meltingtheorem} to the generating function corresponding to their respective RHS.

\section{Preliminaries} \label{sec:preliminaries}

A \emph{partition} of $m\in\Z_{\ge 0}$ is a non-increasing sequence $\lambda=(\lambda_1,\lambda_2,...)$ of non-negative integers summing up to $m$; we denote this by $\lambda\vdash m$. The integer $m$ is called the \emph{size} of $\lambda$ (denoted $|\lambda|$) and the \emph{length} of $\lambda$ (denoted $\ell(\lambda)$) is said to be the number of its non-zero elements.

The \emph{Young diagram with shape} $\lambda$ is the arrangement of $|\lambda|$ $1{\times} 1$ boxes with $\lambda_1$ boxes in the bottom row, $\lambda_2$ in the second, etc., all aligned to the bottom-left (i.e., we use the French notation for the diagrams). Additionally, for a partition $\mu=(\mu_1,\mu_2,...)$ with $\mu_i\le\lambda_i$ for all $i>0$, we define the \emph{skew diagram with shape} $\lambda/\mu$ to be the set difference between the boxes of $\lambda$ and $\mu$.

A \emph{Young tableau} $S$ of shape $\lambda/\mu$ is a filling of the boxes with positive integers. Additionally, if the entries are increasing upwards and non-decreasing rightwards, we call the tableau \emph{semistandard} and denote it with $S\in\SSYT(\lambda/\mu)$. A semistandard Young tableau of size $n$ with entries $1$ through $n$ is called \emph{standard} (denoted $S\in\SYT(\lambda/\mu)$).

Let $\lla/\mmu = (\lambda^1/\mu^1,...,\lambda^m/\mu^m)$ be a sequence of skew partitions\footnote{Note how from now on, we use bolded notation for sequences of shapes and regular font for single shapes.}. We say that a cell $\square=(x,y)$ (i.e., the box in row $y$ and column $x$) of $\lambda^i/\mu^i$ has \emph{content} $c(\square) \coloneqq x-y$ and \emph{shifted content} $\tilde{c}(\square) \coloneqq mc(\square) + i$.

Let $\bm{S}\in\SSYT(\lla/\mmu)$, where $\SSYT(\lla/\mmu)$ denotes the set of semistandard Young tableaux on $\lla/\mmu$. The set of \emph{inversions} in $\bm{S}$ is defined to be $\Inv(bm{S}) \coloneqq \{(\square,\square')\in\lla/\mmu \mid 0 < \tilde{c}(\square')-\tilde{c}(\square) < m \text{ and } \bm{S} (\square)>\bm{S} (\square')\}$. We denote $\inv(\bm{S})\coloneqq |\Inv(\bm{S})|$.

The \textit{LLT polynomial} of $\lla/\mmu$ is the generating function
\begin{equation}
\LLT(\lla/\mmu) \coloneqq \sum_{\bm{S}\in\SSYT(\lla/\mmu)} q^{\inv(\bm{S})}x^{\bm{S}}.
\end{equation}
If $\lla/\mmu$ is a sequence of unicellular shapes, then $\LLT(\lla/\mmu)$ is called a \emph{unicelullar LLT polynomial}.

In \cite{Dolega2019}, Dołęga introduced the notion of Macdonald cumulants: expressions that aim to shed more light on the structure and combinatorics behind Macdonald polynomials. In particular, a proof of Conjecture 1.2 from \cite{Dolega2019} about Schur positivity of Macdonald cumulants would generalize the celebrated result of Haiman \cite{Haiman2001}.

In \cite{DolegaKowalski2021}, the author together with Dołęga, apply similar reasoning to LLT polynomials.

\begin{definition}[For a more general definition, see Definition 2.2 in \cite{DolegaKowalski2021}]
	Let $\lla/\mmu = (\lambda^1/\mu^1,\dots,\lambda^m/\mu^m)$ be a sequence of skew partitions. We define the \emph{LLT cumulant} of $\lla/\mmu$ to be the expression
	\begin{equation} \label{def:cum}
	\ka(\lla/\mmu) \coloneqq (q-1)^{-(m - 1)} \sum\limits_{\mathcal{B}\in\Part(m)} \prod\limits_{B\in\mathcal{B}} (-1)^{|B|-1}(|B|-1)!\LLT(B),
	\end{equation}
	where $\Part(m)$ denotes the set of all set partitions of $[m] \coloneqq \{1,\dots,m\}$, $B$-s are blocks of a set partition $\mathcal{B}$, and 
	\[\LLT(B) \coloneqq \LLT(\lambda^{i_1}/\mu^{i_1},\dots,\lambda^{i_r}/\mu^{i_r})\ \  \text{for}\ B=\{i_1 <\cdots < i_r\}.\]
\end{definition}

Equivalently, using the M\"obius inversion formula for the poset of set partitions, we can define unicellular LLT cumulants recursively by
\begin{equation} \label{def:cumrecursive}
\LLT(\lla/\mmu) = \sum\limits_{\mathcal{B}\in\Part(m)} (q-1)^{m - |\mathcal{B}|} \prod\limits_{B\in\mathcal{B}} \ka(B),
\end{equation}
where
\[\ka(B) \coloneqq \ka(\lambda^{i_1}/\mu^{i_1},\dots,\lambda^{i_r}/\mu^{i_r}) \ \  \text{for}\ B=\{i_1 <\cdots < i_r\}.\]

The inspiration for introducing cumulants to the theory of symmetric functions comes from probability theory: there, the notion allows the study of random variable dependence. In the case of LLT cumulants, their usefulness is best portrayed via so called \emph{LLT graphs}.

\section{LLT graphs} \label{sec:lltgraphs}

LLT graphs were introduced by Dołęga and the author in \cite{DolegaKowalski2021} as a means to understand $\LLT(\lla/\mmu)$ and $\ka(\lla/\mmu)$ as generating functions of combinatorial objects: the graphs themselves.

\begin{definition}
	A directed graph $G$ is called an \emph{LLT graph} if it has three types of edges, visually depicted as $\rightarrow$, $\twoheadrightarrow$, and $\Rightarrow$, and called \emph{edges of type I}, \emph{of type II}, and \emph{double edges}, respectively. Denote the corresponding sets of edges by $E_1(G)$, $E_2(G)$, and $E_d(G)$. 
	
	\emph{The LLT function corresponding to} $G$ is the expression
	\begin{equation} \label{def:graphllt}
	\LLT(G)\coloneqq \sum\limits_{f:V(G)\rightarrow\mathbb{Z}_{>0}} \left(\prod\limits_{(u,v)\in E(G)} \varphi_f(u,v)\right)\cdot \left(\prod\limits_{v\in V(G)} x_{f(v)}\right),
	\end{equation}
	with
	\begin{equation} \label{def:graphcol}
	\varphi_f(u,v)=\begin{cases} [f(u)>f(v)] & \text{for } (u,v)\in E_1(G); \\
	[f(u)\ge f(v)] & \text{for } (u,v)\in E_2(G); \\
	q[f(u)>f(v)]+[f(u)\le f(v)] & \text{for } (u,v)\in E_d(G), \end{cases}
	\end{equation}
	where $[A]$ is the characteristic function of condition $A$, i.e., is equal to $1$ if $A$ is true and $0$ otherwise.
\end{definition}

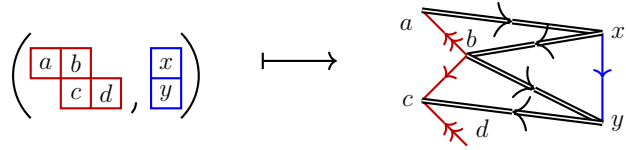
\begin{figure}
	\centering
	\begin{tikzpicture}[scale=.4, every node/.style={scale=0.8}]
	\draw[thick, red] (0,0) -- (-1,0) -- (-1,-1) -- (2,-1) -- (2,-2) -- (0,-2) -- (0,0) -- (1,0) -- (1,-2);
	\draw[thick, blue] (3,-1) -- (3,0) -- (4,0) -- (4,-2) -- (3,-2) -- (3,-1) -- (4,-1);
	\node at (-1,-3) {};
	\draw[thick] (2.5,-1.8) arc (0:-45:.5);
	\draw[thick] (-1,.5) arc (135:225:2);
	\draw[thick] (4,.5) arc (45:-45:2);
	
	\node at (-.5,-.5) {$a$};
	\node at (.5,-.5) {$b$};
	\node at (.5,-1.5) {$c$};
	\node at (1.5,-1.5) {$d$};
	\node at (3.5,-.5) {$x$};
	\node at (3.5,-1.5) {$y$};
	\end{tikzpicture}
	\begin{tikzpicture}[scale=1]
	\draw[thick, |->] (0,0) -- (1,0);
	\node at (-.5,-1) {};
	\end{tikzpicture}
	\hspace{.5cm}
	\begin{tikzpicture}[scale=0.6, every node/.style={scale=0.8}]
	\draw[thick, red] (0,0) -- (.5,-.5);
	\draw[thick, red, <<-] (.5,-.5) -- (1,-1);
	\draw[thick, red, ->] (1,-1) -- (.5,-1.5);
	\draw[thick, red] (.5,-1.5) -- (0,-2);
	\draw[thick, red] (0,-2) -- (.5,-2.5);
	\draw[thick, red, <<-] (.5,-2.5) -- (1,-3);
	
	\draw[thick, blue, ->] (4,-.5) -- (4,-1.5);
	\draw[thick, blue] (4,-1.5) -- (4,-2.5);
	
	\draw[thick, double, ->] (0,0) -- (2,-.25);
	\draw[thick, double] (2,-.25) -- (4,-.5);
	\draw[thick, double, ->] (4,-.5) -- (2.5,-.75);
	\draw[thick, double] (2.5,-.75) -- (1,-1);
	\draw[thick, double, ->] (1,-1) -- (2.5,-1.75);
	\draw[thick, double] (2.5,-1.75) -- (4,-2.5);
	\draw[thick, double, ->] (4,-2.5) -- (2,-2.25);
	\draw[thick, double] (2,-2.25) -- (0,-2);

	\node[below left] at (0,0) {$a$};
	\node[above] at (1.1,-1) {$b$};
	\node[left] at (0,-2) {$c$};
	\node[above right] at (1,-3) {$d$};
	\node[right] at (4,-.5) {$x$};
	\node[right] at (4,-2.5) {$y$};
	\end{tikzpicture}
	\caption{The LLT graph corresponding to ((3,2)/(1), (1,1)).}
	\label{fig:lltgraph}
\end{figure}

Given a sequence of skew shapes $\lla/\mmu$, there is an obvious way to associate with it an LLT graph $G_{\lla/\mmu}$ so that $\LLT(G_{\lla/\mmu}) = \LLT(\lla/\mmu)$. To be precise, we let vertices correspond to cells, edges of type I go from a cell (i.e., the vertex associated with that cell) to the cell directly below it, edges of type II go from a cell to that directly to its left, and double edges connect cells that contribute inversions (see \cref{fig:lltgraph}).

The definitions of the $\LLT$ map and of $\varphi_f$ give simple relations that will prove helpful farther down (note that the below lemma was also proven in \cite{DolegaKowalski2021} using different language).

\begin{lemma} \label{lem:arrowrel}
	\begin{enumerate}
		\item For $i\in\{1,2,d\}$, denote by $G_i$ the LLT graph with $V(G_i) = \{u,v\}$ and $E(G_i) = E_i(G_i) = \{(u,v)\}$. Additionally, let $G_0$ be the graph with $V(G_0) = \{u,v\}$ and $E(G_0) = \emptyset$. Then
		\begin{enumerate}
			\item $\LLT(G_0) = \LLT(G_1) + \LLT(G_2)$,
			\item $\LLT(G_d) = q\LLT(G_1) + \LLT(G_2)$,
			\item $\LLT(G_d) = (q-1)\LLT(G_1) + \LLT(G_0)$.
		\end{enumerate}
		\item For $i\in[2]$, denote by $G_i$ and $G_i'$ the LLT graphs with $V(G_i) = V(G_i') = \{u,v,w\}$, $E(G_i) = E_i(G_i) = \{(u,v),(u,w),(v,w)\}$, and $E(G_i') = E_i(G_i') = \{(u,v),(v,w)\}$. Then $\LLT(G_i) = \LLT(G_i')$, $i=1,2$.
		\item If an LLT graph $G$ contains a directed cycle consisting of edges of type I and type II with at least one edge of type I, then $\LLT(G)=0$.
	\end{enumerate}
\end{lemma}

\begin{proof}
	\begin{enumerate}
		\item
		\begin{enumerate}
			\item By definition, $\LLT(G_0)$ is the sum over all colorings $f$ of two vertices with no restrictions on the colors. We can either have $f(u) > f(v)$ or $f(u) \le f(v)$, which corresponds to $(u,v)\in E_1(G)$ and $(v,u)\in E_2(G)$, respectively.
			\item Similarly to point (a), the coloring $f$ can give either $f(u) > f(v)$ or $f(u) \le f(v)$ (i.e., $(u,v)\in E_1(G)$ or $(v,u)\in E_2(G)$, respectively) with the first contributing $q$ to the summand by definition of double edges and $\varphi_f$.
			\item Follows from (a) and (b).
		\end{enumerate}
		\item Take $i = 2$ (the proof for $i = 1$ is analogous). A coloring $f$ of $G_2$ must satisfy $f(u) \ge f(v)$ and $f(v) \ge f(w)$, which forces $f(u) \ge f(w)$. That, in turn, is equivalent to having $(u,w) \in E_2(G)$.
		\item Let $v_1,...,v_n$ be the subsequent vertices of a directed cycle, i.e., $(v_i,v_{i+1}) \in E_1(G) \cup E_2(G)$ for $i=1,...,n-1$ and $(v_n,v_1)\in E_1(G) \cup E_2(G)$. Then a legal coloring $f$ of $G$ must satisfy $f(v_i) \ge f(v_{i+1})$ for $i=1,...,n-1$ and $f(v_n) \ge f(v_1)$, which can only happen if $f(v_1) = \cdots = f(v_n)$. That, in turn, is impossible if one of the edges is of type I.
	\end{enumerate}
\end{proof}

Note that due to the multiplicative character of \eqref{def:graphllt}, we can apply any of the above relations to an arbitrary LLT graph locally as long as the rest of the graph remains unchanged.

\subsection{Diagrams corresponding to spanning trees} \label{sec:cumunicell}

Recall that \cref{conj:Tdecgen} associates LLT polynomials with trees. We will now describe this correspondence.

Let $T$ be a plane rooted tree. Heuristically, we use an algorithm similar to depth-first search to construct consecutive shapes of $\nu(T)$ by deleting the maximal left-most paths from the root. To be precise, we decompose $V(T) = W_1\cup\cdots\cup W_k$, where $W_i$-s are paths chosen according to the following rules.

Let $W_1$ be the geodesic path from the root to the left-most leaf of $T$. Suppose that we have already defined $W_1,\dots,W_l$. To define $W_{l+1}$, choose a vertex $v\in V(T)\setminus(W_1\cup\cdots\cup W_l)$ following the depth-first search algorithm with respect to $W_1,\dots,W_l$. To be precise, $v$ is the vertex which has a neighbour in some $W_j$ with $1\le j\le l$ maximal and $v$ is the left-most vertex with the maximal distance from the root for the $j$ determined above. We let $W_{l+1}$ be the left-most geodesic path from $v$ downwards to a leaf of $T$.

We will now show that the decomposition is well defined. Firstly, to see that $V(T) = W_1\cup\cdots\cup W_k$, suppose for contradiction that there exists $v\in V(T)\setminus(W_1\cup\cdots\cup W_k)$. Without loss of generality, assume that $v$ is the left-most vertex with that property. Then, $v$ cannot be the root (since the root is in $W_1$), so there exist vertices $v_1,\dots,v_l$ with $v_l=v$, $v_i$ the parent of $v_{i+1}$, $i=1,\dots,l-1$, and $v_i$-s not in $W_1\cup\cdots\cup W_k$. Next, if all vertices to the left of $v_1$ are contained in $W_1\cup\cdots\cup W_r$ for some $1\le r\le k$, then, by the construction above, $v_1$ must be the first vertex of $W_{r+1}$.

Lastly, to see that the decomposition is unique, suppose that we have two decompositions $V(T) = W_1\cup\cdots\cup W_k = W'_1\cup\cdots\cup W'_r$ and that $W_i=W'_i$ for $i=1,\dots,l$. Then, by construction, the first vertex $v$ of $W_{l+1}$ is a neighbour of some $W_j$ with $1\le j\le l$ maximal and similarly for the first vertex $v'$ of $W'_{l+1}$. What is more, both $v$ and $v'$ are the left-most vertices with the maximal distance from the root with those properties, which means that $v=v'$. Since there exists a unique left-most geodesic path from $v$ to a leaf of $T$, we must have $W_{l+1}=W'_{l+1}$.

\begin{definition} \label{def:treetoverticals}
	 With notation as above, denote by $\nu(T)=(\nu^1,\dots,\nu^k)$ the \emph{LLT polynomial associated with the tree} $T$, where $\nu^i=(1^{s_i})/(1^{t_i})$, $i=1,\dots,k$ with
	 \begin{enumerate}
	 	\item $s_i = h(T) - d(W_i)$, where $h(T)=\max\{d(r,v)\mid v\text{ is a leaf in } T\}$ and $d(W)=\min\{d(r,v)\mid v\in W\}$ with $r$ the root of $T$ and $d(r,v)$ the vertex length of the path from $r$ to $v$, and
	 	\item $t_i = s_i - |W_i|$.
	 \end{enumerate}
\end{definition}

\begin{remark}
	Note that the sequence $\nu(T)$ is indeed well defined. Firstly, we always have $s_i > 0$ simply by the definition of $h(T)$. Secondly, $t_i\ge 0$ since the path from the root to the first vertex of $W_i$ followed by $W_i$ is shorter or equal to the maximal path in $T$, and thus $h(T) \ge |W_i|+d(W_i)$, $i=1,\dots,k$.
\end{remark}

Recall that the statements of \cref{thm:meltinglollipopgeneralization} and \cref{thm:singcelldecomp} mention spanning trees of certain graphs with a set of integers as the vertex set rather than plane rooted trees. Fortunately, there is an obvious way to associate with such a graph $T$ a planar drawing in the feel of \cref{def:treetoverticals}. To be precise, we root $T$ in the vertex labeled by the smallest integer and draw the children of each vertex from left to right in an increasing order (see \cref{fig:shapeoftree}).

For simplicity of notation, let us write $\nu(T)$ for the sequence of vertical strips associated to the planar drawing of a tree $T$ with vertex set $[m]$ in the sense of \cref{def:treetoverticals} (see \cref{fig:shapeoftree}). Clearly, such a map $\nu$ is not injective since there exist multiple labelings of a given plane rooted tree which satisfy the rules above. What is more, $\nu$ does not map onto the set of all sequences of vertical-strip shapes since, with notation as in \cref{def:treetoverticals}, $s_1 > s_i$ for $i>1$. In fact, the image of $\nu$ has a particularly nice description in terms of Schr\"oder paths and parking functions, which we now introduce.

\section{Schr\"{o}der paths and parking functions} \label{sec:schroderpaths}

\begin{definition}
	\textit{A Schr\"{o}der path} of length $n$ is a lattice path from $(0,0)$ to $(m,m)$ with steps $n = (0,1)$, $e = (1,0)$, and $d = (1,1)$ (referred to as \textit{north}, \textit{east}, and \textit{diagonal steps}, respectively), which never falls below the main diagonal that joins the ends and has no $d$ steps on that diagonal. We denote by $(i,j)$ the coordinates of the $1\times 1$ box with top right vertex in $(i,j)$.
\end{definition}

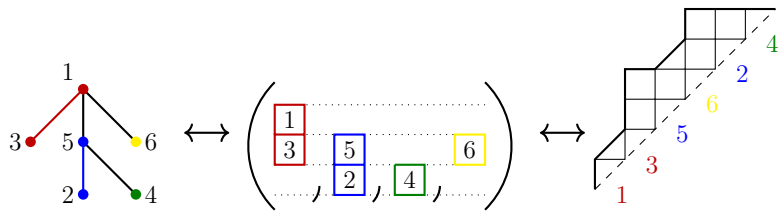
\begin{figure}
	\centering
	\begin{tikzpicture}[scale=.7, every node/.style={scale=0.8}]
	\draw[thick, red] (0,0) -- (-1,-1);
	\draw[thick] (0,0) -- (0,-1);
	\draw[thick] (0,0) -- (1,-1);
	\draw[thick, blue] (0,-1) -- (0,-2);
	\draw[thick] (0,-1) -- (1,-2);
	\fill[red] (0,0) circle (0.1);
	\fill[blue] (0,-1) circle (0.1);
	\fill[blue] (0,-2) circle (0.1);
	\fill[red] (-1,-1) circle (0.1);
	\fill[yellow] (1,-1) circle (0.1);
	\fill[green] (1,-2) circle (0.1);
	\node[above left] at (0,0) {$1$};
	\node[left] at (0,-1) {$5$};
	\node[left] at (0,-2) {$2$};
	\node[left] at (-1,-1) {$3$};
	\node[right] at (1,-1) {$6$};
	\node[right] at (1,-2) {$4$};
	\end{tikzpicture}
	\begin{tikzpicture}[scale=1.2]
	\draw[thick, <->] (0,0) -- (.5,0);
	\node at (0,-.75) {};
	\end{tikzpicture}
	\begin{tikzpicture}[scale=.6, every node/.style={scale=0.8}]
	\draw[dotted] (-.66,0) -- (4,0);
	\draw[dotted] (-.66,-.66) -- (4,-.66);
	\draw[dotted] (-.66,-1.33) -- (4,-1.33);
	\draw[dotted] (-.66,-2) -- (4,-2);
	\draw[thick, red] (-.66,-.66) -- (-.66,0) -- (0,0) -- (0,-1.33) -- (-.66,-1.33) -- (-.66,-.66) -- (0,-.66);
	\draw[thick] (.33,-1.8) arc (0:-45:.5);
	\draw[thick, blue] (.66,-1.33) -- (.66,-.66) -- (1.33,-.66) -- (1.33,-2) -- (.66,-2) -- (.66,-1.33) -- (1.33,-1.33);
	\draw[thick] (1.66,-1.8) arc (0:-45:.5);
	\draw[thick, green] (2,-1.33) -- (2.66,-1.33) -- (2.66,-2) -- (2,-2) -- (2,-1.33);
	\draw[thick] (3,-1.8) arc (0:-45:.5);
	\draw[thick, yellow] (3.33,-1.33) -- (4,-1.33) -- (4,-.66) -- (3.33,-.66) -- (3.33,-1.33);
	\draw[thick] (-.66,.5) arc (135:225:2);
	\draw[thick] (4,.5) arc (45:-45:2);
	\node at (-.33,-.33) {$1$};
	\node at (-.33,-1) {$3$};
	\node at (1,-1) {$5$};
	\node at (1,-1.66) {$2$};
	\node at (2.33,-1.66) {$4$};
	\node at (3.66,-1) {$6$};
	\end{tikzpicture}
	\begin{tikzpicture}[scale=1.2]
	\draw[thick, <->] (0,0) -- (.5,0);
	\node at (0,-.75) {};
	\end{tikzpicture}
	\begin{tikzpicture}[scale=.8, every node/.style={scale=0.8}]
	\draw[dashed] (0,0) -- (3,3);
	\draw (0,.5) -- (.5,.5) -- (.5,1) -- (1,1) -- (1,1.5) -- (1.5,1.5) -- (1.5,2) -- (2,2) -- (2,2.5) -- (2.5,2.5) -- (2.5,3);
	\draw (.5,1.5) -- (1,1.5) -- (1,2) -- (1.5,2) -- (1.5,2.5) -- (2,2.5) -- (2,3);
	\draw[thick] (0,0) -- (0,.5) -- (.5,1) -- (.5,2) -- (1,2) -- (1.5,2.5) -- (1.5,3) -- (3,3);
	\node[red, below right] at (.2,.2) {$1$};
	\node[red, below right] at (.7,.7) {$3$};
	\node[blue, below right] at (1.2,1.2) {$5$};
	\node[yellow, below right] at (1.7,1.7) {$6$};
	\node[blue, below right] at (2.2,2.2) {$2$};
	\node[green, below right] at (2.7,2.7) {$4$};
	\end{tikzpicture}
	\caption{A spanning tree with its corresponding sequence of shapes and Schr\"{o}der path.}
	\label{fig:shapeoftree}
\end{figure}

For such a path $P$, we say that its $i$-th column \emph{is of height} $h(i)=k$ if the point $(i-1,k)$ lies on $P$ and is followed by either an east or diagonal step, $1\le i\le m$. We define the \emph{jump} of the $i$-th column to be the value $j(i)\coloneqq h(i+1)-h(i)$, $i=1,\dots,m-1$. Lastly, a box $(i,j)$ is called \textit{an outer corner} of $P$ if the point $(i-1,j-1)$ lies on $P$ and is followed by a sequence of an east step and a north step.

\begin{proposition} \label{prop:verticalstripstoschroderpaths}
	For every $m\in\mathbb{Z}_{>0}$, there exists a bijection between the set of Schr\"oder paths of length $m$ and the set of LLT polynomials of sequences of vertical-strip shapes with $m$ boxes.
\end{proposition}

The above result is a classical correspondence between Schr\"oder paths and LLT polynomials of vertical-strip shapes (see, e.g., \cite{Haglund2007}). Below, we present an alternative proof using LLT graphs.

\begin{proof}[Proof of \cref{prop:verticalstripstoschroderpaths}]
	We will construct a bijection between Schr\"oder paths of length $m$ and certain LLT graphs, which we will later identify with LLT polynomials of sequences of vertical-strip shapes according to \cref{sec:lltgraphs}.
	
	For a path $P$ of length $m$, label box $(i,i)$ by the integer $i$, $1\le i\le m$. Then, construct the LLT graph $G(P)$ with $V(G(P))=[m]$ and $E_2(G(P))=\emptyset$ according to the following rules:
	\begin{enumerate}
		\item if $P$ has a diagonal step in the box $(i,j)$, then $(i,j)\in E_1(G(P))$;
		\item if the box $(i,j)$, $i<j$, lies under $P$, then $(i,j)\in E_d(G(P))$.
	\end{enumerate}
	
	It is easy to see that the map $P\longmapsto G(P)$ is injective. Indeed, if $P$ and $P'$ agree on the first $k$ steps and differ in the $(k+1)$-th step which begins at point $(i,j)$, then either
	\begin{enumerate}
		\item $(i,j)\notin E(G(P))$ and ($(i,j)\in E_1(G(P'))$ or $(i,j)\in E_d(G(P'))$),
		\item $(i,j)\notin E(G(P'))$ and ($(i,j)\in E_1(G(P))$ or $(i,j)\in E_d(G(P))$),
		\item $(i,j)\in E_1(G(P))$ and $(i,j)\in E_d(G(P'))$,
		\item $(i,j)\in E_d(G(P))$ and $(i,j)\in E_1(G(P'))$.
	\end{enumerate}
	
	Observe that $G(P)$ restricted to only edges of type I is a set of disjoint directed paths, and thus corresponds to vertical strips in the sense of \cref{sec:lltgraphs}. Furthermore, $E_d(G(P))$ is in a one-to-one correspondence with inversion pairs between vertical strips. Indeed, for two type I directed paths $I=(i_1,\dots,i_s)$ and $J=(j_1,\dots,j_t)$ with $i_1<j_1$, we must have a directed path $K=(k_1,\dots,k_r)$ of double edges alternating between consecutive vertices of $I$ and $J$ with $k_1=\max\{i_k\mid k\in [s], i_k< j_1\}$ and either
	\[k_r=\begin{cases}
		\min\{i_k\mid k\in[s], i_k>j_t\} & \text{for } i_s>j_t, \\
		\min\{j_k\mid k\in[t], j_k>i_s\} & \text{otherwise.}
	\end{cases}\]
	
	To obtain the inverse map, it is enough to reverse the procedure.
\end{proof}

Let us label the bijection in \cref{prop:verticalstripstoschroderpaths} by $\mu:P\longmapsto \mu(P)$. It turns out that we can connect $\mu$ with the map $\nu$ we introduced in \cref{sec:cumunicell}.

\begin{proposition} \label{prop:treestoschroder}
	There exists a bijection $T\leftrightarrow P$ between the set of plane rooted trees with $m$ vertices and the set of Schr\"oder paths of length $m$ which satisfy the following conditions:
	\begin{enumerate}
		\item $P$ is connected (i.e., the only points of $P$ on the main diagonal are $(0,0)$ and $(m,m)$);
		\item the first steps of $P$ are $n$ and $d$; and
		\item $P$ has no outer corners,
	\end{enumerate}
	such that $\nu(T)=\mu(P)$.
\end{proposition}
\begin{proof}
	We claim that the composition $\mu^{-1}\circ\nu$ is the bijection in question. Indeed, $\mu^{-1}\circ\nu$ is an injection from planar rooted trees with $n$ vertices to Schr\"oder paths of length $n$. Therefore, it is enough to prove that the composition maps onto the set of Schr\"oder paths that satisfy the conditions from the statement.
	
	Let $T$ be a plane rooted tree. The bijection given in \cref{def:treetoverticals} gives a sequence of vertical-strip shapes $\nu(T)$. Let $P$ be the Schr\"oder path that satisfies $\nu(T)=\mu(P)$.
	
	The connectedness of $T$ translates straightforwardly to condition $(1)$ of $P$. Also, observe that the map $T\longmapsto \nu(T)$ forces the cell corresponding to the root $r\in T$ to be the unique cell $\nu(r)$ with the smallest content in $\nu(T)$, and thus, it corresponds to box $(1,1)$ in $\mu(P)$. Furthermore, $r$'s left-most child maps to the cell below $\nu(r)$ in $(\nu(T))^1$, which translates to $P$ satisfying condition $(2)$ from the statement.
	
	To see that $P$ satisfies condition $(3)$, suppose for contradiction that $P$ has an outer corner in the box $(i,j)$ with no outer corners in boxes $(x,y)$ with $x<i$. This means that in the LLT graph $G=G(P)$ corresponding to $\mu(P)$, we have $(i,i+1),(j-1,j),(i,j-1),(i+1,j-1),(i+1,j)\in E_d(G)$ but $(i,j)\notin E(G)$. Visually, this translates to the subgraph
	
	
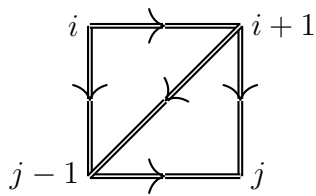
\begin{center}
	\begin{tikzpicture}
		\draw[thick, double, ->] (4,-.5) -- (4,-1.5);
		\draw[thick, double] (4,-1.5) -- (4,-2.5);
		\draw[thick, double, ->] (2,-.5) -- (2,-1.5);
		\draw[thick, double] (2,-1.5) -- (2,-2.5);
		
		\draw[thick, double, ->] (2,-.5) -- (3,-.5);
		\draw[thick, double] (3,-.5) -- (4,-.5);
		\draw[thick, double, ->] (2,-2.5) -- (3,-2.5);
		\draw[thick, double] (3,-2.5) -- (4,-2.5);

		\draw[thick, double, ->] (4,-.5) -- (3,-1.5);
		\draw[thick, double] (3,-1.5) -- (2,-2.5);
		
		\node[left] at (2,-.5) {$i$};
		\node[left] at (2,-2.5) {$j-1$};
		\node[right] at (4,-.5) {$i+1$};
		\node[right] at (4,-2.5) {$j$};
	\end{tikzpicture}
	\end{center}
	which further translates to the four possible cell positionings in $\mu(P)$ (here, cells at the same level have the same content in $\mu(P)$):
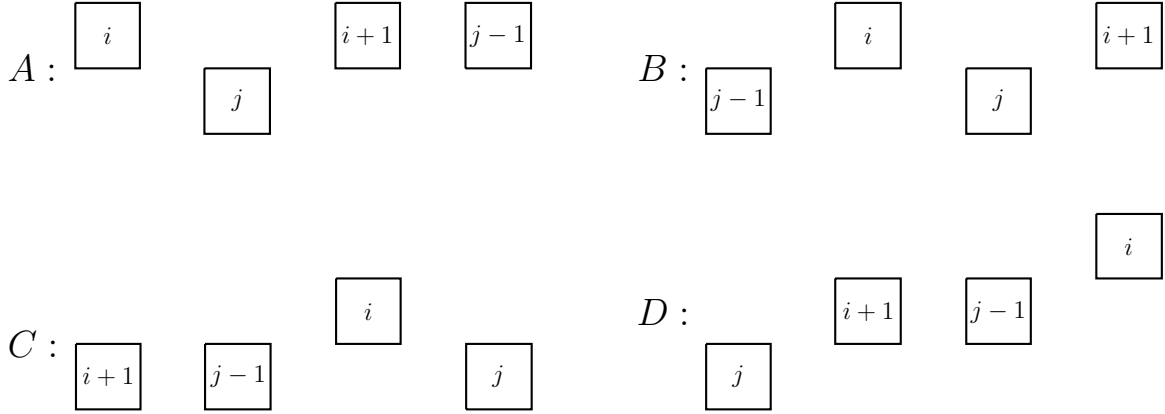
	\begin{center}
	\begin{tikzpicture}[scale=1.3, every node/.style={scale=0.8}]
	\draw[thick] (-.66,-.66) -- (0,-.66) -- (0,-1.33) -- (-.66,-1.33) -- (-.66,-.66);
	\draw[thick] (.66,-1.33) -- (1.33,-1.33) -- (1.33,-2) -- (.66,-2) -- (.66,-1.33);
	\draw[thick] (2,-.66) -- (2.66,-.66) -- (2.66,-1.33) -- (2,-1.33) -- (2,-.66);
	\draw[thick] (3.33,-1.33) -- (4,-1.33) -- (4,-.66) -- (3.33,-.66) -- (3.33,-1.33);
	\node[scale=1.5] at (-1.1,-1.33) {$A:$};
	\node at (-.33,-1) {$i$};
	\node at (1,-1.66) {$j$};
	\node at (2.33,-1) {$i+1$};
	\node at (3.66,-1) {$j-1$};
	\end{tikzpicture}
	\hspace{1cm}
	\begin{tikzpicture}[scale=1.3, every node/.style={scale=0.8}]
	\draw[thick] (-.66,-1.33) -- (0,-1.33) -- (0,-2) -- (-.66,-2) -- (-.66,-1.33);
	\draw[thick] (.66,-.66) -- (1.33,-.66) -- (1.33,-1.33) -- (.66,-1.33) -- (.66,-.66);
	\draw[thick] (2,-1.33) -- (2.66,-1.33) -- (2.66,-2) -- (2,-2) -- (2,-1.33);
	\draw[thick] (3.33,-.66) -- (4,-.66) -- (4,-1.33) -- (3.33,-1.33) -- (3.33,-.66);
	\node[scale=1.5] at (-1.1,-1.33) {$B:$};
	\node at (-.33,-1.66) {$j-1$};
	\node at (1,-1) {$i$};
	\node at (2.33,-1.66) {$j$};
	\node at (3.66,-1) {$i+1$};
	\end{tikzpicture}
	
	\vspace{1cm}
	
	\begin{tikzpicture}[scale=1.3, every node/.style={scale=0.8}]
	\draw[thick] (-.66,-1.33) -- (0,-1.33) -- (0,-2) -- (-.66,-2) -- (-.66,-1.33);
	\draw[thick] (.66,-1.33) -- (1.33,-1.33) -- (1.33,-2) -- (.66,-2) -- (.66,-1.33);
	\draw[thick] (2,-.66) -- (2.66,-.66) -- (2.66,-1.33) -- (2,-1.33) -- (2,-.66);
	\draw[thick] (3.33,-2) -- (4,-2) -- (4,-1.33) -- (3.33,-1.33) -- (3.33,-2);
	\node[scale=1.5] at (-1.1,-1.33) {$C:$};
	\node at (-.33,-1.66) {$i+1$};
	\node at (1,-1.66) {$j-1$};
	\node at (2.33,-1) {$i$};
	\node at (3.66,-1.66) {$j$};
	\end{tikzpicture}
	\hspace{1cm}
	\begin{tikzpicture}[scale=1.3, every node/.style={scale=0.8}]
	\draw[thick] (-.66,-1.33) -- (0,-1.33) -- (0,-2) -- (-.66,-2) -- (-.66,-1.33);
	\draw[thick] (.66,-.66) -- (1.33,-.66) -- (1.33,-1.33) -- (.66,-1.33) -- (.66,-.66);
	\draw[thick] (2,-.66) -- (2.66,-.66) -- (2.66,-1.33) -- (2,-1.33) -- (2,-.66);
	\draw[thick] (3.33,0) -- (4,-0) -- (4,-.66) -- (3.33,-.66) -- (3.33,0);
	\node[scale=1.5] at (-1.1,-1) {$D:$};
	\node at (-.33,-1.66) {$j$};
	\node at (1,-1) {$i+1$};
	\node at (2.33,-1) {$j-1$};
	\node at (3.66,-.33) {$i$};
	\end{tikzpicture}
	\end{center}
	In particular, there exists no vertex $v\in G(P)$ such that $(i,v)\in E_1(G)$ or $(v,j)\in E_1(G)$.
	
	In cases $A$, $B$, and $C$, the position of $j$ means that in $T$, $j$ is a child of $i$. But by the definition of $\nu$, this contradicts the absence of a vertex $v\in G(P)$ such that $(i,v)\in E_1(G)$. In the case of $D$, the parent of $j$ in $T$ would have content between $c(i)$ and $c(j)$, which would place the vertex between $i$ and $i+1$, which is impossible.
	
	To see that $\mu^{-1}\circ\nu:T\longmapsto P$ is surjective, it is enough to reverse the reasoning above. I.e., if $P$ is a Schr\"oder path satisfying conditions $(1)$-$(3)$ from the statement, the conditions ensure that $\mu(P)=\nu(T)$ for some tree $T$.
\end{proof}

For an example of the correspondence between Cayley trees and Schr\"oder paths satisfying the relations, see the first two arrows \cref{fig:shapeoftree}.

As a matter of fact, \cref{prop:treestoschroder} gives the possibility to express the LLT polynomials corresponding to special sequences of vertical-strip shapes using another class of well-known combinatorial objects.

A \emph{Dyck path} is a Schr\"oder path with no diagonal steps. Observe that for every Schr\"oder path $P$, we can easily associate with it a Dyck path $D(P)$ by exchanging each diagonal step of $P$ by a sequence of an east step and a north step. Furthermore, if $P$ satisfies conditions from the statement of \cref{prop:treestoschroder}, the first three steps of $D(P)$ are determined: a north step, an east step, and a north step. Thus, we obtain a simple bijection between Schr\"oder paths of length $m$ that satisfy conditions $(1)$-$(3)$ from \cref{prop:treestoschroder} and Dyck paths of length $m-1$.

\subsection{Schr\"oder path relations}

Similarly to how we introduced relations between graphs (or rather between the generating functions of the graphs), we can study relations between Schr\"{o}der paths. This approach was used by Alexandersson and Sulzgruber \cite{AlexanderssonSulzgruber2020}, who show the following result (also known as "local linear relations" in \cite{Lee2018}).

\begin{lemma} \label{lem:schroderrelations}
	Let $P$ be a Schr\"{o}der path.
	\begin{enumerate}[label=(\Alph*)]
		\item If $P=SneT$ for some paths $S$ and $T$, then
		\[\LLT(SneT) = (q-1)\LLT(SdT) + \LLT(SenT).\]
		\item If $P = SndReeT$ for some paths $S$, $T$, and $R$ with $S$ ending in $(i,j)$ and $SndR$ ending in $(j,k)$, then
		\[\LLT(SndReeT) = q\LLT(SdnReeT).\]
	\end{enumerate}
\end{lemma}

The visual representation of the above relations in the spirit of \cref{lem:arrowrel} takes the form
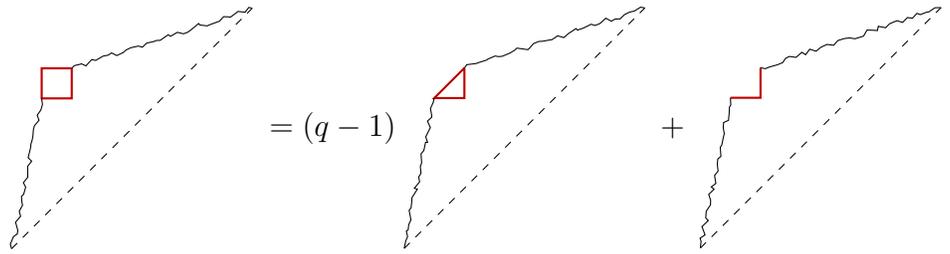
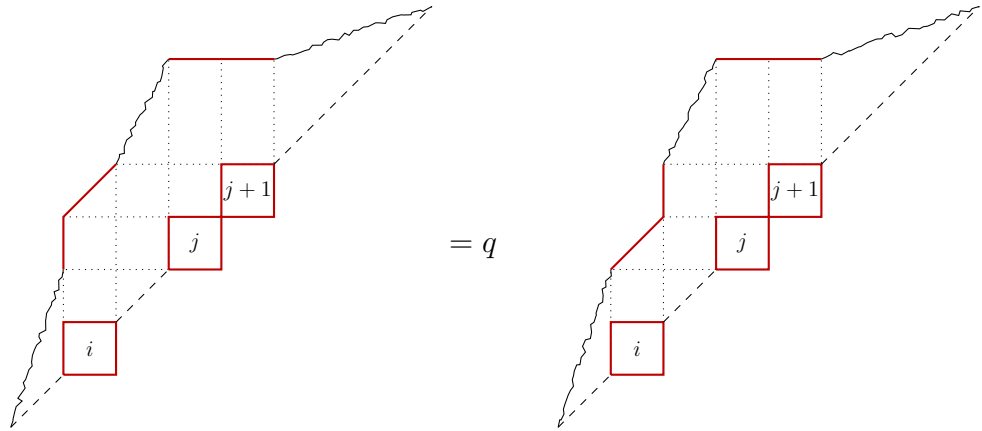
\begin{enumerate}
	\item
	\begin{equation*}
	\vcenter{
		\hbox{
			\begin{tikzpicture}[scale=0.8]
			\draw[thick, red] (1,3) -- (.5,3) -- (.5,2.5) -- (1,2.5) -- (1,3);
			\draw[dashed] (0,0) -- (4,4);
			\path[draw,clip,decoration={random steps,segment length=2pt,amplitude=1pt}] decorate {(1,3) -- (4,4)} decorate {(0,0) -- (0.5,2.5)};
			\end{tikzpicture}
		}
	}
	=
	(q-1)\vcenter{
		\hbox{
			\begin{tikzpicture}[scale=0.8]
			\draw[thick, red] (1,3) -- (.5,2.5) -- (1,2.5) -- (1,3);
			\draw[dashed] (0,0) -- (4,4);
			\path[draw,clip,decoration={random steps,segment length=2pt,amplitude=1pt}] decorate {(1,3) -- (4,4)} decorate {(0,0) -- (0.5,2.5)};
			\end{tikzpicture}
		}
	}
	+
	\vcenter{
		\hbox{
			\begin{tikzpicture}[scale=0.8]
			\draw[thick, red] (1,3) -- (1,2.5) -- (.5,2.5);
			\draw[dashed] (0,0) -- (4,4);
			\path[draw,clip,decoration={random steps,segment length=2pt,amplitude=1pt}] decorate {(1,3) -- (4,4)} decorate {(0,0) -- (0.5,2.5)};
			\end{tikzpicture}
		}
	},
	\end{equation*}
	
	\item
	\begin{equation*}
	\vcenter{
		\hbox{
			\begin{tikzpicture}[scale=1.4, every node/.style={scale=0.7}]
			\draw[dotted] (1,2.5) -- (2.5,2.5) -- (2.5,4);
			\draw[dotted] (1,1) -- (1,2) -- (2,2) -- (2,4);
			\draw[dotted] (1.5,1.5) -- (1.5,3) -- (3,3) -- (3,4);
			\draw[thick, red] (1,2) -- (1,2.5) -- (1.5,3);
			\draw[thick, red] (2,4) -- (3,4);
			\draw[thick, red] (1,1) -- (1,1.5) -- (1.5,1.5) -- (1.5,1) -- (1,1);
			\draw[thick, red] (2,2) -- (2,2.5) -- (3,2.5) -- (3,3) -- (2.5,3) -- (2.5,2) -- (2,2);
			\draw[dashed] (.5,.5) -- (1,1);
			\draw[dashed] (1.5,1.5) -- (2,2);
			\draw[dashed] (3,3) -- (4.5,4.5);
			\node at (1.25,1.25) {$i$};
			\node at (2.25,2.25) {$j$};
			\node at (2.75,2.75) {$j+1$};
			\path[draw,clip,decoration={random steps,segment length=2pt,amplitude=1pt}] decorate {(3,4) -- (4.5,4.5)} decorate {(.5,.5) -- (1,2)} decorate {(1.5,3) -- (2,4)};
			\end{tikzpicture}
		}
	}
	=
	q\vcenter{
		\hbox{
			\begin{tikzpicture}[scale=1.4, every node/.style={scale=0.7}]
			\draw[dotted] (1.5,2.5) -- (2.5,2.5) -- (2.5,4);
			\draw[dotted] (1,1) -- (1,2) -- (2,2) -- (2,4);
			\draw[dotted] (1.5,1.5) -- (1.5,3) -- (3,3) -- (3,4);
			\draw[thick, red] (1,2) -- (1.5,2.5) -- (1.5,3);
			\draw[thick, red] (2,4) -- (3,4);
			\draw[thick, red] (1,1) -- (1,1.5) -- (1.5,1.5) -- (1.5,1) -- (1,1);
			\draw[thick, red] (2,2) -- (2,2.5) -- (3,2.5) -- (3,3) -- (2.5,3) -- (2.5,2) -- (2,2);
			\draw[dashed] (.5,.5) -- (1,1);
			\draw[dashed] (1.5,1.5) -- (2,2);
			\draw[dashed] (3,3) -- (4.5,4.5);
			\node at (1.25,1.25) {$i$};
			\node at (2.25,2.25) {$j$};
			\node at (2.75,2.75) {$j+1$};
			\path[draw,clip,decoration={random steps,segment length=2pt,amplitude=1pt}] decorate {(3,4) -- (4.5,4.5)} decorate {(.5,.5) -- (1,2)} decorate {(1.5,3) -- (2,4)};
			\end{tikzpicture}
		}
	}.
	\end{equation*}
\end{enumerate}

Due to the correspondence between Schr\"oder paths and LLT graphs, these identities translate to identities on the corresponding LLT graphs. To simplify the notation, whenever we state a relation of the form $\sum_{i=1}^n  f_i(q)\cdot G_i = 0$, where $f_i\in \mathbb{Z}[q]$ and $G_i$ are LLT graphs, we understand it as $\sum_{i=1}^n  f_i(q)\cdot \LLT(G_i) = 0$. What is more, to simplify the formulas even further, we depict pictorially only the subgraph of the bigger graph that involves some changes on the edges. For instance, the second identity locally (i.e., limited to the three vertices associated with the marked boxes) has the following pictorial form:
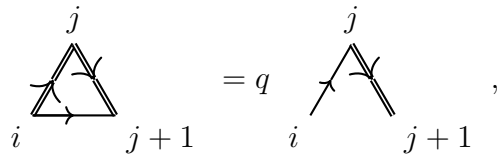
\begin{equation*}
\vcenter{
	\hbox{
		\begin{tikzpicture}[scale=1.1]
		\draw[thick,double, ->] (-.5,0) -- (-.25,.43);
		\draw[thick,double] (-.25,.43) -- (0,.86);
		\draw[thick,double, ->] (0,.86) -- (.25,.43);
		\draw[thick,double] (.25,.43) -- (.5,0);
		\draw[thick, ->] (-0.5,0) -- (0,0);
		\draw[thick] (-0.5,0) -- (0.5,0);
		\node[below left] at (-0.5,0) {$i$};
		\node[above] at (0,0.86) {$j$};
		\node[below right] at (0.5,0) {$j+1$};
		\end{tikzpicture}
	}
}
=q
\vcenter{
	\hbox{
		\begin{tikzpicture}[scale=1.1]
		\draw[thick, ->] (-.5,0) -- (-.25,.43);
		\draw[thick] (-.25,.43) -- (0,.86);
		\draw[thick,double, ->] (0,.86) -- (.25,.43);
		\draw[thick,double] (.25,.43) -- (.5,0);
		\node[below left] at (-0.5,0) {$i$};
		\node[above] at (0,0.86) {$j$};
		\node[below right] at (0.5,0) {$j+1$};
		\end{tikzpicture}
	}
},
\end{equation*}
which formally means that for an LLT graph $G$ associated with the Schr\"oder path $P$ from point $(B)$ of \cref{lem:schroderrelations}, one has
\[\LLT(G) = q\LLT(G'),\]
where the only differences between $G$ and $G'$ are on the subgraph induced by the vertex set $\{i,j,j+1\}$ and are described by the picture.

We present a different proof of \cref{lem:schroderrelations} from the one in \cite{AlexanderssonSulzgruber2020}. It is an easy manipulation of the relations presented in \cref{lem:arrowrel} on graphs corresponding to Schr\"oder paths satisfying the statement of \cref{lem:schroderrelations} and thus, the relations hold true for those special cases of LLT graphs.

\begin{proof}
	\begin{enumerate}
		\item The identity is a straightforward translation of property $(C)$ from \cref{lem:arrowrel}.
		
		\item From the shape of $P$ we deduce that the only relations that differ between the left- and right-hand sides are between cells $i$, $j$, and $j+1$. Therefore, we first focus on the objects locally and later extend to the whole graphs.
		
		If we use graph relations from \cref{lem:arrowrel}, we get (for simplicity of notation, we drop the vertex labels)
		\begin{align*}
		\vcenter{\hbox{
				\begin{tikzpicture}[scale=1.1]
				\draw[thick,double, ->] (-.5,0) -- (-.25,.43);
				\draw[thick,double] (-.25,.43) -- (0,.86);
				\draw[thick,double, ->] (0,.86) -- (.25,.43);
				\draw[thick,double] (.25,.43) -- (.5,0);
				\draw[thick, ->] (-0.5,0) -- (0,0);
				\draw[thick] (-0.5,0) -- (0.5,0);
				\end{tikzpicture}}}
		&\hspace{.4cm}\stackrel{(1c)}{=}\hspace{.4cm}
		(q-1)\vcenter{\hbox{
				\begin{tikzpicture}[scale=1.1]
				\draw[thick, ->] (-.5,0) -- (-.25,.43);
				\draw[thick] (-.25,.43) -- (0,.86);
				\draw[thick,double, ->] (0,.86) -- (.25,.43);
				\draw[thick,double] (.25,.43) -- (.5,0);
				\draw[thick, ->] (-0.5,0) -- (0,0);
				\draw[thick] (-0.5,0) -- (0.5,0);
				\end{tikzpicture}}}
		+
		\vcenter{\hbox{
				\begin{tikzpicture}[scale=1.1]
				\draw[thick,double, ->] (0,.86) -- (.25,.43);
				\draw[thick,double] (.25,.43) -- (.5,0);
				\draw[thick, ->] (-0.5,0) -- (0,0);
				\draw[thick] (-0.5,0) -- (0.5,0);
				\end{tikzpicture}}}
		\\
		&\hspace{.4cm}\stackrel{(1a)}{=}\hspace{.4cm}
		(q-1)\vcenter{\hbox{
				\begin{tikzpicture}[scale=1.1]
				\draw[thick, ->] (-.5,0) -- (-.25,.43);
				\draw[thick] (-.25,.43) -- (0,.86);
				\draw[thick,double, ->] (0,.86) -- (.25,.43);
				\draw[thick,double] (.25,.43) -- (.5,0);
				\end{tikzpicture}}}
		-
		(q-1)\vcenter{\hbox{
				\begin{tikzpicture}[scale=1.1]
				\draw[thick, ->] (-.5,0) -- (-.25,.43);
				\draw[thick] (-.25,.43) -- (0,.86);
				\draw[thick,double, ->] (0,.86) -- (.25,.43);
				\draw[thick,double] (.25,.43) -- (.5,0);
				\draw[thick] (-0.5,0) -- (-.1,0);
				\draw[thick, <<-] (-0.1,0) -- (0.5,0);
				\end{tikzpicture}}}
		+
		\vcenter{\hbox{
				\begin{tikzpicture}[scale=1.1]
				\draw[thick,double, ->] (0,.86) -- (.25,.43);
				\draw[thick,double] (.25,.43) -- (.5,0);
				\filldraw[black] (-0.5,0) circle (1pt);
				\end{tikzpicture}}}
		-
		\vcenter{\hbox{
				\begin{tikzpicture}[scale=1.1]
				\draw[thick,double, ->] (0,.86) -- (.25,.43);
				\draw[thick,double] (.25,.43) -- (.5,0);
				\draw[thick] (-0.5,0) -- (-.1,0);
				\draw[thick, <<-] (-0.1,0) -- (0.5,0);
				\end{tikzpicture}}}
		\hspace{1cm} \\
		&\stackrel{(1c)+(1a)}{=}
		(q-1)\vcenter{\hbox{
				\begin{tikzpicture}[scale=1.1]
				\draw[thick, ->] (-.5,0) -- (-.25,.43);
				\draw[thick] (-.25,.43) -- (0,.86);
				\draw[thick,double, ->] (0,.86) -- (.25,.43);
				\draw[thick,double] (.25,.43) -- (.5,0);
				\end{tikzpicture}}}
		-
		(q-1)\vcenter{\hbox{
				\begin{tikzpicture}[scale=1.1]
				\draw[thick, ->] (-.5,0) -- (-.25,.43);
				\draw[thick] (-.25,.43) -- (0,.86);
				\draw[thick, ->] (0,.86) -- (.25,.43);
				\draw[thick] (.25,.43) -- (.5,0);
				\draw[thick] (-0.5,0) -- (-.1,0);
				\draw[thick, <<-] (-0.1,0) -- (0.5,0);
				\end{tikzpicture}}}
		-
		(q-1)\vcenter{\hbox{
				\begin{tikzpicture}[scale=1.1]
				\draw[thick, ->] (-.5,0) -- (-.25,.43);
				\draw[thick] (-.25,.43) -- (0,.86);
				\draw[thick] (-0.5,0) -- (-.1,0);
				\draw[thick, <<-] (-0.1,0) -- (0.5,0);
				\end{tikzpicture}}}
		\\
		&\hspace{1.4cm}+\hspace{.4cm}
		\vcenter{\hbox{
				\begin{tikzpicture}[scale=1.1]
				\draw[thick, ->] (-.5,0) -- (-.25,.43);
				\draw[thick] (-.25,.43) -- (0,.86);
				\draw[thick,double, ->] (0,.86) -- (.25,.43);
				\draw[thick,double] (.25,.43) -- (.5,0);
				\end{tikzpicture}}}
		+
		\vcenter{\hbox{
				\begin{tikzpicture}[scale=1.1]
				\draw[thick] (-.5,0) -- (-.3,.336);
				\draw[thick, <<-] (-.3,.336) -- (0,.86);
				\draw[thick,double, ->] (0,.86) -- (.25,.43);
				\draw[thick,double] (.25,.43) -- (.5,0);
				\end{tikzpicture}}}
		-
		\vcenter{\hbox{
				\begin{tikzpicture}[scale=1.1]
				\draw[thick,double, ->] (0,.86) -- (.25,.43);
				\draw[thick,double] (.25,.43) -- (.5,0);
				\draw[thick] (-0.5,0) -- (-.1,0);
				\draw[thick, <<-] (-0.1,0) -- (0.5,0);
				\end{tikzpicture}}}
		\\
		&\stackrel{(3)+(1c)}{=}
		q\vcenter{\hbox{
				\begin{tikzpicture}[scale=1.1]
				\draw[thick, ->] (-.5,0) -- (-.25,.43);
				\draw[thick] (-.25,.43) -- (0,.86);
				\draw[thick,double, ->] (0,.86) -- (.25,.43);
				\draw[thick,double] (.25,.43) -- (.5,0);
				\end{tikzpicture}}}
		-
		(q-1)\vcenter{\hbox{
				\begin{tikzpicture}[scale=1.1]
				\draw[thick, ->] (-.5,0) -- (-.25,.43);
				\draw[thick] (-.25,.43) -- (0,.86);
				\draw[thick] (-0.5,0) -- (-.1,0);
				\draw[thick, <<-] (-0.1,0) -- (0.5,0);
				\end{tikzpicture}}}
		\\
		&\hspace{1.4cm}+\hspace{.4cm}
		(q-1)\vcenter{\hbox{
				\begin{tikzpicture}[scale=1.1]
				\draw[thick] (-.5,0) -- (-.3,.336);
				\draw[thick, <<-] (-.3,.336) -- (0,.86);
				\draw[thick, ->] (0,.86) -- (.25,.43);
				\draw[thick] (.25,.43) -- (.5,0);
				\end{tikzpicture}}}
		+
		\vcenter{\hbox{
				\begin{tikzpicture}[scale=1.1]
				\draw[thick] (-.5,0) -- (-.3,.336);
				\draw[thick, <<-] (-.3,.336) -- (0,.86);
				\filldraw[black] (0.5,0) circle (1pt);
				\end{tikzpicture}}}
		-
		(q-1)\vcenter{\hbox{
				\begin{tikzpicture}[scale=1.1]
				\draw[thick, ->] (0,.86) -- (.25,.43);
				\draw[thick] (.25,.43) -- (.5,0);
				\draw[thick] (-0.5,0) -- (-.1,0);
				\draw[thick, <<-] (-0.1,0) -- (0.5,0);
				\end{tikzpicture}}}
		-
		\vcenter{\hbox{
				\begin{tikzpicture}[scale=1.1]
				\draw[thick] (-0.5,0) -- (-.1,0);
				\draw[thick, <<-] (-0.1,0) -- (0.5,0);
				\filldraw[black] (0,0.86) circle (1pt);
				\end{tikzpicture}}}
		\\
		&\hspace{.4cm}\stackrel{(1a)}{=}\hspace{.4cm}
		q\vcenter{\hbox{
				\begin{tikzpicture}[scale=1.1]
				\draw[thick, ->] (-.5,0) -- (-.25,.43);
				\draw[thick] (-.25,.43) -- (0,.86);
				\draw[thick,double, ->] (0,.86) -- (.25,.43);
				\draw[thick,double] (.25,.43) -- (.5,0);
				\end{tikzpicture}}}
		-
		(q-1)\vcenter{\hbox{
				\begin{tikzpicture}[scale=1.1]
				\draw[thick, ->] (-.5,0) -- (-.25,.43);
				\draw[thick] (-.25,.43) -- (0,.86);
				\draw[thick] (-0.5,0) -- (-.1,0);
				\draw[thick, <<-] (-0.1,0) -- (0.5,0);
				\end{tikzpicture}}}
		+
		(q-1)\vcenter{\hbox{
				\begin{tikzpicture}[scale=1.1]
				\draw[thick] (-.5,0) -- (-.3,.336);
				\draw[thick, <<-] (-.3,.336) -- (0,.86);
				\draw[thick, ->] (0,.86) -- (.25,.43);
				\draw[thick] (.25,.43) -- (.5,0);
				\draw[thick, ->] (-0.5,0) -- (0,0);
				\draw[thick] (-0.5,0) -- (0.5,0);
				\end{tikzpicture}}}
		+
		(q-1)\vcenter{\hbox{
				\begin{tikzpicture}[scale=1.1]
				\draw[thick] (-.5,0) -- (-.3,.336);
				\draw[thick, <<-] (-.3,.336) -- (0,.86);
				\draw[thick, ->] (0,.86) -- (.25,.43);
				\draw[thick] (.25,.43) -- (.5,0);
				\draw[thick] (-0.5,0) -- (-.1,0);
				\draw[thick, <<-] (-0.1,0) -- (0.5,0);
				\end{tikzpicture}}}
		\\
		&\hspace{1.4cm}+\hspace{.4cm}
		\vcenter{\hbox{
				\begin{tikzpicture}[scale=1.1]
				\draw[thick] (-.5,0) -- (-.3,.336);
				\draw[thick, <<-] (-.3,.336) -- (0,.86);
				\filldraw[black] (0.5,0) circle (1pt);
				\end{tikzpicture}}}
		-
		(q-1)\vcenter{\hbox{
				\begin{tikzpicture}[scale=1.1]
				\draw[thick, ->] (0,.86) -- (.25,.43);
				\draw[thick] (.25,.43) -- (.5,0);
				\draw[thick] (-0.5,0) -- (-.1,0);
				\draw[thick, <<-] (-0.1,0) -- (0.5,0);
				\end{tikzpicture}}}
		-
		\vcenter{\hbox{
				\begin{tikzpicture}[scale=1.1]
				\draw[thick] (-0.5,0) -- (-.1,0);
				\draw[thick, <<-] (-0.1,0) -- (0.5,0);
				\filldraw[black] (0,0.86) circle (1pt);
				\end{tikzpicture}}}
		\\
		&\hspace{.4cm}=\hspace{.4cm}
		q\vcenter{\hbox{
				\begin{tikzpicture}[scale=1.1]
				\draw[thick, ->] (-.5,0) -- (-.25,.43);
				\draw[thick] (-.25,.43) -- (0,.86);
				\draw[thick,double, ->] (0,.86) -- (.25,.43);
				\draw[thick,double] (.25,.43) -- (.5,0);
				\end{tikzpicture}}}.
		\end{align*}
		The last equality follows from
		\begin{equation*}
		\vcenter{\hbox{
				\begin{tikzpicture}[scale=1.1]
				\draw[thick] (-.5,0) -- (-.3,.336);
				\draw[thick, <<-] (-.3,.336) -- (0,.86);
				\draw[thick, ->] (0,.86) -- (.25,.43);
				\draw[thick] (.25,.43) -- (.5,0);
				\draw[thick] (-0.5,0) -- (-.1,0);
				\draw[thick, <<-] (-0.1,0) -- (0.5,0);
				\end{tikzpicture}}}
		-
		\vcenter{\hbox{
				\begin{tikzpicture}[scale=1.1]
				\draw[thick, ->] (0,.86) -- (.25,.43);
				\draw[thick] (.25,.43) -- (.5,0);
				\draw[thick] (-0.5,0) -- (-.1,0);
				\draw[thick, <<-] (-0.1,0) -- (0.5,0);
				\end{tikzpicture}}}
		=
		0
		\end{equation*}
		(which is a result of \cref{lem:arrowrel}) and
		\begin{equation*}
		\vcenter{\hbox{
				\begin{tikzpicture}[scale=1.1]
				\draw[thick] (-.5,0) -- (-.3,.336);
				\draw[thick, <<-] (-.3,.336) -- (0,.86);
				\filldraw[black] (0.5,0) circle (1pt);
				\end{tikzpicture}}}
		-
		\vcenter{\hbox{
				\begin{tikzpicture}[scale=1.1]
				\draw[thick] (-0.5,0) -- (-.1,0);
				\draw[thick, <<-] (-0.1,0) -- (0.5,0);
				\filldraw[black] (0,0.86) circle (1pt);
				\end{tikzpicture}}}
		=
		0,
		\hspace{2cm}
		\vcenter{\hbox{
				\begin{tikzpicture}[scale=1.1]
				\draw[thick, ->] (-.5,0) -- (-.25,.43);
				\draw[thick] (-.25,.43) -- (0,.86);
				\draw[thick] (-0.5,0) -- (-.1,0);
				\draw[thick, <<-] (-0.1,0) -- (0.5,0);
				\end{tikzpicture}}}
		-
		\vcenter{\hbox{
				\begin{tikzpicture}[scale=1.1]
				\draw[thick] (-.5,0) -- (-.3,.336);
				\draw[thick, <<-] (-.3,.336) -- (0,.86);
				\draw[thick, ->] (0,.86) -- (.25,.43);
				\draw[thick] (.25,.43) -- (.5,0);
				\draw[thick, ->] (-0.5,0) -- (0,0);
				\draw[thick] (-0.5,0) -- (0.5,0);
				\end{tikzpicture}}}
		=
		0
		\end{equation*}
		which is more subtle and results from the following two properties of our graph.
		\begin{enumerate}
			\item Denote the four graphs from left to right by $G_1$, $G_2$, $G_3$, and $G_4$. Then, for every legal coloring $f$ (i.e., having non-zero contribution $\varphi_f$ into \cref{def:graphllt}, therefore satisfying comparability relations implied by \cref{def:graphcol}) of $G_1$ (respectively, $G_3$), we can define a legal coloring $\tilde{f}$ of $G_2$ (respectively, $G_4$) by setting $\tilde{f}(i) = f(i)$, $\tilde{f}(j) = f(j+1)$, and $\tilde{f}(j+1) = f(j)$. What is more, one easily verifies that the map $f\to \tilde{f}$ is a bijection preserving the contribution $\varphi$ into \cref{def:graphllt}.
			\item For every $v>j+1$, the pair $(j,v)$ is an edge of type $A\in\{\rightarrow,\twoheadrightarrow,\Rightarrow\}$ if and only if $(j+1,v)$ is an edge of the same type $A$.
		\end{enumerate}
	
		Therefore, we obtain a bijection between the colorings $f$ of the whole graph on $[n]$ given in the statement with $G|_{\{i,j,j+1\}} = G_1$ (respectively, $G|_{\{i,j,j+1\}} = G_3$) and the colorings $\tilde{f}$ of the same graph with $G|_{\{i,j,j+1\}} = G_2$ (respectively, $G|_{\{i,j,j+1\}} = G_4$) by setting $\tilde{f}(k)$ as in (a) above for $k \in \{i,j,j+1\}$ and $\tilde{f}(k) = f(k)$ otherwise. Moreover, this bijection preserves the weight $\varphi$ and the associated colors and thus, $\LLT(G_1) = \LLT(G_2)$ and $\LLT(G_3) = \LLT(G_4)$.
	\end{enumerate}
\end{proof}

\section{Melting lollipops and proof of \cref{thm:meltinglollipopgeneralization}} \label{sec:meltinglollipops}

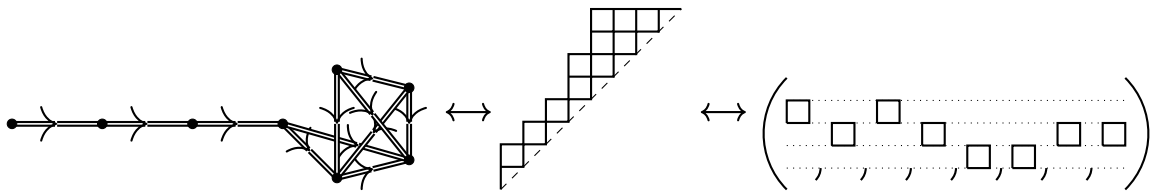
\begin{figure}
	\centering
	\begin{tikzpicture}[scale=1.2]
			\draw[double,thick, ->] (0,0) -- (.5,0);
			\draw[double,thick] (.5,0) -- (1,0);
			\draw[double,thick, ->] (1,0) -- (1.5,0);
			\draw[double,thick] (1.5,0) -- (2,0);
			\draw[double,thick, ->] (2,0) -- (2.5,0);
			\draw[double,thick] (2.5,0) -- (3,0);
			\filldraw[black] (0,0) circle (1.5pt);
			\filldraw[black] (1,0) circle (1.5pt);
			\filldraw[black] (2,0) circle (1.5pt);
			
			\draw[double,thick, ->] (3,0) -- (3.3,-.3);
			\draw[double,thick] (3.3,-.3) -- (3.6,-.6);
			\draw[double,thick, ->] (3.6,-.6) -- (4,-.5);
			\draw[double,thick] (4,-.5) -- (4.4,-.4);
			\draw[double,thick, ->] (3,0) -- (3.84,-.24);
			\draw[double,thick] (3.84,-.24) -- (4.4,-.4);
			\draw[double,thick, ->] (3.6,.6) -- (4,.5);
			\draw[double,thick] (4,.5) -- (4.4,.4);
			\draw[double,thick, ->] (4.4,.4) -- (4.4,0);
			\draw[double,thick] (4.4,0) -- (4.4,-.4);
			\draw[double,thick, ->] (4.4,.4) -- (4,-.1);
			\draw[double,thick] (4,-.1) -- (3.6,-.6);
			\draw[double,thick, ->] (3.6,.6) -- (4,.1);
			\draw[double,thick] (4,.1) -- (4.4,-.4);
			\draw[double,thick, ->] (3.6,.6) -- (3.6,0);
			\draw[double,thick] (3.6,0) -- (3.6,-.6);
			\filldraw[black] (3,0) circle (1.5pt);
			\filldraw[black] (3.6,-.6) circle (1.5pt);
			\filldraw[black] (3.6,.6) circle (1.5pt);
			\filldraw[black] (4.4,.4) circle (1.5pt);
			\filldraw[black] (4.4,-.4) circle (1.5pt);
	\end{tikzpicture}
	\begin{tikzpicture}[scale=1.2]
		\draw[thick, <->] (0,0) -- (.5,0);
		\node at (0,-.75) {};
	\end{tikzpicture}
	\begin{tikzpicture}[scale=.3]
			\draw[thick] (0,0) -- (0,2) -- (1,2) -- (1,3) -- (2,3) -- (2,4) -- (3,4) -- (3,6) -- (4,6) -- (4,8) -- (8,8);
			\draw[thick] (0,1) -- (1,1) -- (1,2) -- (2,2) -- (2,3) -- (3,3) -- (3,4) -- (4,4) -- (4,5) -- (5,5) -- (5,6) -- (6,6) -- (6,7) -- (7,7) -- (7,8);
			\draw[thick] (3,5) -- (4,5) -- (4,6) -- (5,6) -- (5,7) -- (6,7) -- (6,8);
			\draw[thick] (4,7) -- (5,7) -- (5,8);
			\draw[dashed] (0,0) -- (8,8);
	\end{tikzpicture}
	\begin{tikzpicture}[scale=1.2]
		\draw[thick, <->] (0,0) -- (.5,0);
		\node at (0,-.75) {};
	\end{tikzpicture}
	\begin{tikzpicture}[scale=.3]
		\draw[thick] (0,2) -- (0,3) -- (1,3) -- (1,2) -- (0,2);
		\draw[thick] (2,1) -- (2,2) -- (3,2) -- (3,1) -- (2,1);
		\draw[thick] (4,2) -- (4,3) -- (5,3) -- (5,2) -- (4,2);
		\draw[thick] (6,1) -- (6,2) -- (7,2) -- (7,1) -- (6,1);
		\draw[thick] (8,0) -- (8,1) -- (9,1) -- (9,0) -- (8,0);
		\draw[thick] (10,0) -- (10,1) -- (11,1) -- (11,0) -- (10,0);
		\draw[thick] (12,1) -- (12,2) -- (13,2) -- (13,1) -- (12,1);
		\draw[thick] (14,1) -- (14,2) -- (15,2) -- (15,1) -- (14,1);
		\draw[thick] (15,4) arc (45:-45:3.5);
		\draw[thick] (0,4) arc (135:225:3.5);
		\draw[thick] (1.5,0) arc (0:-45:.75);
		\draw[thick] (3.5,0) arc (0:-45:.75);
		\draw[thick] (5.5,0) arc (0:-45:.75);
		\draw[thick] (7.5,0) arc (0:-45:.75);
		\draw[thick] (9.5,0) arc (0:-45:.75);
		\draw[thick] (11.5,0) arc (0:-45:.75);
		\draw[thick] (13.5,0) arc (0:-45:.75);
		\draw[dotted] (0,0) -- (15,0);
		\draw[dotted] (0,1) -- (15,1);
		\draw[dotted] (0,2) -- (15,2);
		\draw[dotted] (0,3) -- (15,3);
	\end{tikzpicture}
	\caption{The melting lollipop (i.e., its LLT graph variant) $L^{(2)}_{(3,5)}$ and its corresponding Dyck path and sequence of shapes.}
	\label{fig:meltinglollipopgraph}
\end{figure}

Recall from \cref{def:meltinglollipop} that we obtain the melting lollipop $L^{(k)}_{(l,m)}$ by joining $K_m$ to a path of length $l$ and erasing $k$ edges (see \cref{fig:meltinglollipopgraph} for an example). In fact, there is a straightforward way to translate the notion to the language of LLT functions. To be precise, we will denote by $L^{(k)}_{(l,m)}$ the LLT graph on the vertex set $[l+m]$ with $E_1(L^{(k)}_{(l,m)}) = E_2(L^{(k)}_{(l,m)}) = \emptyset$ and $E_d(L^{(k)}_{(l,m)})$ equal to the set $\{(i,i+1) \mid i\in[l]\} \cup (\{(i,j)\mid i,j\in [l+1,l+m], i<j\} \setminus \{(l+1,l+m),\dots(l+1,l+m-k+1)\})$.

In what follows, we only use the LLT variants of melting lollipops and thus, we use for them the same notation as in \cref{def:meltinglollipop}. What is more, it is easy to see that melting lollipops correspond to a class of LLT polynomials of unicellular diagrams and those, in turn, correspond to a class of Dyck paths via the map described in \cref{sec:lltgraphs} (see \cref{fig:meltinglollipopgraph} for an example of the correspondence).

Before we move on to the proof of \cref{thm:meltinglollipopgeneralization}, let us introduce some notation. Let $\lla/\mmu$ be a sequence of $m$ unicellular shapes and denote the cells by $\square_1,\dots,\square_m$ with $\tilde{c}(\square_1) < \dots < \tilde{c}(\square_m)$. For a subset $B\subseteq [m]$, write $G_{\lla/\mmu}(B)$ for the non-directed graph with $V(G_{\lla/\mmu}(B))=B$ and $\{i,j\}\in E(G_{\lla/\mmu}(B))$ whenever $0 < |\tilde{c}(\square_i) - \tilde{c}(\square_j)| < m$.

\subsection{Proof of \cref{thm:meltinglollipopgeneralization}} \label{sec:proof}

\begin{proof}
	Using \cref{def:cumrecursive}, we have
	\[\LLT(\lla/\mmu) = \sum_{\mathcal{B}\in\Part(l+m)} (q-1)^{l+m-|\mathcal{B}|} \prod_{B\in\mathcal{B}} \kappa(B),\]
	where $\kappa(B) = \kappa(\la^{i_1}/\mu^{i_1},\dots,\la^{i_r}/\mu^{i_r} \mid B = \{i_1<\dots<i_r\})$.
	
	However, as shown in \cite[Theorem 3.8]{DolegaKowalski2021}, $\kappa(B) = 0$ whenever $G_{\lla/\mmu}(B)$ is disconnected. Thus, \eqref{eq:meltingtheorem} is, in fact, equivalent to
	\[\begin{aligned} \LLT(\lla/\mmu) &= \sum_{\substack{\mathcal{B}\in\Part(l+m) \\ B\in\mathcal{B} \Rightarrow G_{\lla/\mmu}(B) \text{ conncected}}} (q-1)^{l+m-|\mathcal{B}|} \prod_{B\in\mathcal{B}} \kappa(B) \\
	&= \sum_{\substack{\mathcal{B}\in\Part(l+m) \\ B\in\mathcal{B} \Rightarrow G_{\lla/\mmu}(B) \text{ conncected}}} (q-1)^{l+m-|\mathcal{B}|} \prod_{B\in\mathcal{B}} \sum_{T\subseteq K_B} \LLT(\nu(T)),\end{aligned}\]
	where the inner sum runs over all spanning trees of $G_{\lla/\mmu}(B)$ and we construct $\nu(T)$ according to \cref{def:treetoverticals}.
	
	Lastly, it suffices to interpret the product above as $\LLT(\nu(F))$ with $F$ the forest corresponding to the trees on the components $B\in\mathcal{B}$ to get
	\begin{equation} \label{thm:formulawithforests}
	\LLT(\lla/\mmu)(q+1) = \sum\limits_F q^{l+m-\#F}\LLT(\nu(F))(q+1),
	\end{equation}
	where the sum runs over all spanning forests of $G_{\lla/\nnu}([l+m])$, $\#F$ denotes the number of connected components of $F$, and $\LLT(\nu(F))$ is the product of $\LLT(\nu(T_i))$ with $T_i$ the connected components of $F$, each rooted in its vertex of the smallest label.
	
	The idea behind the proof is to take the Schr\"oder path corresponding to $\lla/\mmu$ and repeatedly apply \cref{lem:schroderrelations} to decompose it into a sum of Schr\"oder paths corresponding to spanning forests of $G_{\lla/\nnu}([l+m])$. To be precise, we will show that
	\[\LLT(\lla/\mmu)(q+1) = \sum_P \varphi(P;q) \LLT(\mu(P))(q+1),\]
	where $\varphi(P;q)$ are polynomials in $q$ dependent on the path $P$ and the sum runs over Schr\"oder paths corresponding to spanning forests of $G_{\lla/\mmu}([l+m])$ with each connected component rooted. The coefficients of $\varphi(P;q)$ will correspond to the number of labelings of such forests with integers $1$ through $l+m$ which satisfy the conditions from \cref{sec:cumunicell}. To be precise, we will describe the coefficients using the shape of $P$ and polynomials
	\[W_s(q) \coloneqq \sum_{i=1}^{\tilde{h}(s)-h(s)} {\tilde{h}(s)-h(s) \choose i} q^i\in\Z_{>0}[q]\]
	where $l+1\le s \le l+m$ and $\tilde{h}(s)$ and $h(s)$ denote the heights of column $s$ in $L^{(k)}_{(l,m)}$ and $P$, respectively.
	
	First of all, we apply relation $(A)$ from \cref{lem:schroderrelations} to each of the first $l$ columns of $L^{(k)}_{(l,m)}$ to reduce each of those columns to minimal height and each diagonal step that the relation gives contributes the factor of $q$ to its corresponding summand. To be precise, we obtain
	\begin{equation} \label{eq:meltingtheoremintermediate}
		\LLT(\lla/\mmu)(q+1) = \sum_P q^{d(P)}\LLT(\mu(P))(q+1),
	\end{equation}
	where the sum runs over all Schr\"oder paths $P$ with $h(i) = i$ for $i=1,\dots,l$, $h(l+1)=l+m-k$, $h(l+2)=l+m$, and $P$ does not have a diagonal step in the $(l+1)$-th column and $d(P)$ is the number of diagonal steps in $P$.
	
	Let $P$ be a Schr\"oder path appearing in \eqref{eq:meltingtheoremintermediate}. We will now describe a procedure that utilizes relations $(A)$ and $(B)$ from \cref{lem:schroderrelations} to decompose the last $m$ columns of $P$. Since it does not affect the first $l$ columns, in the formulas that follow, we draw only the remaining ones.
	
	We begin by reducing the height of the $(l+1)$-th column of $P$. First, we apply relation $(A)$ from \cref{lem:schroderrelations} to the cell $(l+1,l+m-k)$, which gives
	\begin{equation} \label{thm:reduce1stcolumn}
	\vcenter{
		\hbox{
			\begin{tikzpicture}[scale=.5]
			\draw[thick] (0,0) -- (0,2) -- (0.5,2)-- (0.5,3) -- (3,3);
			\draw (0,1.5) -- (0.5,1.5) -- (0.5,2);
			\draw[dashed] (0,0) -- (3,3);
			\end{tikzpicture}
		}
	}
	=
	q\vcenter{
		\hbox{
			\begin{tikzpicture}[scale=.5]
			\draw[thick] (0,0) -- (0,1.5) -- (0.5,2) -- (0.5,3) -- (3,3);
			\draw (0,1.5) -- (0.5,1.5) -- (0.5,2);
			\draw[dashed] (0,0) -- (3,3);
			\end{tikzpicture}
		}
	}
	+
	\vcenter{
		\hbox{
			\begin{tikzpicture}[scale=.5]
			\draw[thick] (0,0) -- (0,1.5) -- (0.5,1.5) -- (0.5,3) -- (3,3);
			\draw[dashed] (0,0) -- (3,3);
			\end{tikzpicture}
		}
	}.
	\end{equation}
	
	Next, to the first summand, we apply relation $(B)$ from \cref{lem:schroderrelations} to move the diagonal step downwards, i.e., to obtain a new shape with a coefficient in $\mathbb{Z}_{>0}[q]$ and with the diagonal step lower than initially. We repeat the procedure consecutively until the diagonal step is in the cell $(l+1,l+2)$.
	
	To the second summand from \eqref{thm:reduce1stcolumn}, we apply $(A)$ from \cref{lem:schroderrelations} to cell $(l+1,l+m-k-1)$. This will again give us two summands, one with a diagonal step in the first column and one without. To these two, we repeat the reasoning presented until we reduce the first column of the path until no longer possible. Therefore, we end up with
	\begin{equation} \label{thm:firstcolumndone}
	\vcenter{
		\hbox{
			\begin{tikzpicture}[scale=.5]
			\draw[thick] (0,0) -- (0,2) -- (0.5,2)-- (0.5,3) -- (3,3);
			\draw[dashed] (0,0) -- (3,3);
			\end{tikzpicture}
		}
	}
	=
	W_{l+1}(q)\vcenter{
		\hbox{
			\begin{tikzpicture}[scale=.5]
			\draw[thick] (0,0) -- (0,0.5) -- (0.5,1) -- (0.5,3) -- (3,3);
			\draw[dashed] (0,0) -- (3,3);
			\end{tikzpicture}
		}
	}
	+
	\vcenter{
		\hbox{
			\begin{tikzpicture}[scale=.5]
			\draw[thick] (0,0) -- (0,0.5) -- (0.5,0.5) -- (0.5,3) -- (3,3);
			\draw[dashed] (0,0) -- (3,3);
			\end{tikzpicture}
		}
	},
	\end{equation}
	where
	\begin{equation*}
	\begin{aligned} W_{l+1}(q) &= q(q+1)^{\tilde{h}(l+1)-h(l+1)-2} + q(q+1)^{\tilde{h}(l+1)-h(l+1)-3} + \cdots + q \\
	&= (q+1)^{\tilde{h}(l+1)-h(l+1)-1} - 1 \\
	&= \sum_{i=1}^{\tilde{h}(l+1)-h(l+1)} {\tilde{h}(l+1)-h(l+1) \choose i}q^i. \end{aligned}
	\end{equation*}
	
	The next steps of the decomposition differ for every summand. For the second shape in \eqref{thm:firstcolumndone}, we reduce the $(l+2)$-th column as we did above. For the first, the exponent of $q$ will determine the jump of the $(l+1)$-th column in the next steps of the decomposition. To be precise, for $1\le s\le m-k-j(l)$, we decompose the summand with $q^{s}$ until the $(l+2)$-th column is of height $l+s+1$ (equivalently, until $j(l+1)=s$). In other words, we use the two identities from \cref{lem:schroderrelations} to get
	\begin{align*}
	{\tilde{h}(l+1)-h(l+1) \choose j(l+1)}&q^{j(l+1)}\vcenter{
		\hbox{
			\begin{tikzpicture}[scale=.5]
			\draw[thick] (0,0) -- (0,0.5) -- (0.5,1) -- (0.5,3) -- (3,3);
			\draw (0,0) -- (0,0.5) -- (0.5,0.5) -- (0.5,1);
			\draw[dashed] (0,0) -- (3,3);
			\end{tikzpicture}
		}
	}
	= \\
	&=
	{\tilde{h}(l+1)-h(l+1) \choose j(l+1)}q^{j(l+1)}W_{l+2}(q)\vcenter{
		\hbox{
			\begin{tikzpicture}[scale=.5]
			\draw[thick] (0,0) -- (0,0.5) -- (0.5,1) -- (0.5,2) -- (1,2.5) -- (1,3) -- (3,3);
			\draw (0.5,2) -- (1,2) -- (1,2.5);
			\draw[dashed] (0,0) -- (3,3);
			\end{tikzpicture}
		}
	} \\
	&+
	{\tilde{h}(l+1)-h(l+1) \choose j(l+1)}q^{j(l+1)} \vcenter{
		\hbox{
			\begin{tikzpicture}[scale=.5]
			\draw[thick] (0,0) -- (0,0.5) -- (0.5,1) -- (0.5,2) -- (1,2) -- (1,3) -- (3,3);
			\draw[dashed] (0,0) -- (3,3);
			\end{tikzpicture}
		}
	},
	\end{align*}
	where, analogously to the previous step,
	\begin{equation*}
	W_{l+2}(q) = (q+1)^{\tilde{h}(l+2)-h(l+2)} - 1 = \sum_{i=1}^{\tilde{h}(l+2)-h(l+2)} {\tilde{h}(l+2)-h(l+2) \choose i}q^i.
	\end{equation*}
	
	In general, assume that we have decomposed the first $2\le s\le m-1$ columns, i.e., that we have managed to express the left-hand side of \eqref{thm:formulawithforests} as a sum of Schr\"{o}der paths whose first $l+s$ columns have been reduced.
	
	Take a summand with coefficient $q^K$ of a path with diagonal steps in columns $l+1\le i_1 < \dots < i_r = s$, $K=j(i_1)+\cdots+j(i_{r-1})+t$, where $1 \leq t \leq m-h(s)$ (if $i_r < s$ then $K=j(i_1)+\cdots+j(i_r)$ and we reduce the $(s+1)$-th column maximally as in the first step). We lower the $(s+1)$-th column to height $l+K+1$.
	
	\begin{align*}
	&\prod_{b=1}^{r-1} {\tilde{h}(i_b) - h(i_b) \choose j(i_b)}\cdot { \tilde{h}(i_r)-h(i_r) \choose t} q^K\vcenter{
		\hbox{
			\begin{tikzpicture}[scale=.5]
			\draw[thick] (0.5,1.5) -- (1,2) -- (1,4) -- (4,4);
			\draw[dashed] (0,0) -- (4,4);
			\path[draw,clip,decoration={random steps,segment length=2pt,amplitude=1pt}] decorate {(0,0) -- (0.5,1.5)};
			\end{tikzpicture}
		}
	}
	= \\
	&=\prod_{b=1}^r {\tilde{h}(i_b) - h(i_b) \choose j(i_b)} q^KW_{i_r+1}(q)\vcenter{
		\hbox{
			\begin{tikzpicture}[scale=.5]
			\draw[thick] (0.5,1.5) -- (1,2) -- (1,2.5) -- (1.5,3) -- (1.5,4) -- (4,4);
			\draw (1,2.5) -- (1.5,2.5) -- (1.5,3);
			\draw[dashed] (0,0) -- (4,4);
			\path[draw,clip,decoration={random steps,segment length=2pt,amplitude=1pt}] decorate {(0,0) -- (0.5,1.5)};
			\end{tikzpicture}
		}
	}
	+\prod_{b=1}^r {\tilde{h}(i_b) - h(i_b) \choose j(i_b)} q^K\vcenter{
		\hbox{
			\begin{tikzpicture}[scale=.5]
			\draw[thick] (0.5,1.5) -- (1,2) -- (1,2.5) -- (1.5,2.5) -- (1.5,4) -- (4,4);
			\draw[dashed] (0,0) -- (4,4);
			\path[draw,clip,decoration={random steps,segment length=2pt,amplitude=1pt}] decorate {(0,0) -- (0.5,1.5)};
			\end{tikzpicture}
		}
	},
	\end{align*}
	where $K = j(i_1) + \cdots + j(i_r)$ (meaning that $t = j(i_r)$) and
	\[W_{i_r+1}(q) = (q+1)^{\tilde{h}(i_r+1) - h(i_r+1)} - 1 = \sum_{i=1}^{\tilde{h}(i_r+1) - h(i_r+1)} {\tilde{h}(i_r+1) - h(i_r+1) \choose i} q.\]
	
	In the end, the above decomposition gives
	\begin{equation} \label{thm:decompositioninpaths}
	\LLT(\lla/\mmu)(q+1) = \sum_P \varphi(P;q)\LLT(\mu(P))(q+1),
	\end{equation}
	where the sum runs over all Schr\"oder paths $P$ contained in $L^{(k)}_{(l,m)}$ (i.e., ones that never go above $L^{(k)}_{(l,m)}$) with $h(l+1) = l+1$ and whose each connected component of length at least $2$ (i.e., a section of the path from one point on the diagonal to another without any others in between) begins with $nd$ and has no outer corners, and
	\begin{equation*}
	\varphi(P;q) = q^{l+m-c(P)}\prod_{s=1}^{l}{ n-h(i_s) \choose n- h(i_s+1)},
	\end{equation*}
	where $c(P)$ is the number of connected components of $P$ and $i_1,\dots,i_r$ are the columns of $P$ containing diagonal steps.
	
	Therefore, each shape appearing in \eqref{thm:decompositioninpaths} corresponds to a forest in the sense of \cref{sec:schroderpaths}. Furthermore, a path $P$ having diagonal steps in columns $i_1,\dots,i_r$ corresponds to the vertex $i_s$ having $j(i_s)$ children in the forest.
	
	It remains to show that the number of forests on $[l+m]$ with labelings satisfying the conditions from \cref{sec:cumunicell} which correspond to the same planar drawing $F$ is equal to the multinomial coefficient from $\varphi(P;q)$ with $P$ the Schr\"oder path associated with $F$. In fact, it is a straightforward consequence of the combinatorics of set partitions and thus, we leave the details to the reader.
\end{proof}
	
\section{Complete graph case and parking functions}
	
Observe that a special case of a melting lollipop is the complete graph $K_m = L^{(0)}_{(0,m)}$, which corresponds to the full Dyck path $P_m = n^me^m$ and to the sequence $\lla/\mmu = ((1),\dots,(1))$ of $m$ unicellular non-skew diagrams.
	
A straightforward application of \cref{thm:meltinglollipopgeneralization} to $K_m$ implies the first equality in \cref{thm:singcelldecomp}. Indeed, in this case, we can understand spanning trees as Cayley trees using the correspondence from \cref{sec:cumunicell}. What is more, this allows for linking the result to yet another family of well-known combinatorial objects: \emph{parking functions}.
	
\begin{definition}
	For $m\in\Z_{>0}$, we say that $f:[m]\rightarrow[m]$ is a \emph{parking function on $m$ cars} if $|f^{-1}([i])|\ge i$ for $i=1,\dots,m$.
\end{definition}
	
Recall that the maps $\nu$ and $\mu$ give simple ways to translate the vertex labels on a Cayley tree $T$ to cell labels in the sequence $\nu(T)$ and to diagonal box labels in the Schr\"oder path $P$ such that $\nu(T)=\mu(P)$ (see \cref{fig:shapeoftree}).

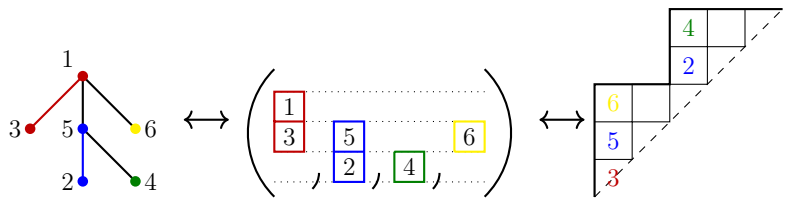
\begin{figure}
	\centering
	\begin{tikzpicture}[scale=.7, every node/.style={scale=0.8}]
	\draw[thick, red] (0,0) -- (-1,-1);
	\draw[thick] (0,0) -- (0,-1);
	\draw[thick] (0,0) -- (1,-1);
	\draw[thick, blue] (0,-1) -- (0,-2);
	\draw[thick] (0,-1) -- (1,-2);
	\fill[red] (0,0) circle (0.1);
	\fill[blue] (0,-1) circle (0.1);
	\fill[blue] (0,-2) circle (0.1);
	\fill[red] (-1,-1) circle (0.1);
	\fill[yellow] (1,-1) circle (0.1);
	\fill[green] (1,-2) circle (0.1);
	\node[above left] at (0,0) {$1$};
	\node[left] at (0,-1) {$5$};
	\node[left] at (0,-2) {$2$};
	\node[left] at (-1,-1) {$3$};
	\node[right] at (1,-1) {$6$};
	\node[right] at (1,-2) {$4$};
	\end{tikzpicture}
	\begin{tikzpicture}[scale=1.2]
	\draw[thick, <->] (0,0) -- (.5,0);
	\node at (0,-.75) {};
	\end{tikzpicture}
	\begin{tikzpicture}[scale=.6, every node/.style={scale=0.8}]
	\draw[dotted] (-.66,0) -- (4,0);
	\draw[dotted] (-.66,-.66) -- (4,-.66);
	\draw[dotted] (-.66,-1.33) -- (4,-1.33);
	\draw[dotted] (-.66,-2) -- (4,-2);
	\draw[thick, red] (-.66,-.66) -- (-.66,0) -- (0,0) -- (0,-1.33) -- (-.66,-1.33) -- (-.66,-.66) -- (0,-.66);
	\draw[thick] (.33,-1.8) arc (0:-45:.5);
	\draw[thick, blue] (.66,-1.33) -- (.66,-.66) -- (1.33,-.66) -- (1.33,-2) -- (.66,-2) -- (.66,-1.33) -- (1.33,-1.33);
	\draw[thick] (1.66,-1.8) arc (0:-45:.5);
	\draw[thick, green] (2,-1.33) -- (2.66,-1.33) -- (2.66,-2) -- (2,-2) -- (2,-1.33);
	\draw[thick] (3,-1.8) arc (0:-45:.5);
	\draw[thick, yellow] (3.33,-1.33) -- (4,-1.33) -- (4,-.66) -- (3.33,-.66) -- (3.33,-1.33);
	\draw[thick] (-.66,.5) arc (135:225:2);
	\draw[thick] (4,.5) arc (45:-45:2);
	\node at (-.33,-.33) {$1$};
	\node at (-.33,-1) {$3$};
	\node at (1,-1) {$5$};
	\node at (1,-1.66) {$2$};
	\node at (2.33,-1.66) {$4$};
	\node at (3.66,-1) {$6$};
	\end{tikzpicture}
	\begin{tikzpicture}[scale=1.2]
	\draw[thick, <->] (0,0) -- (.5,0);
	\node at (0,-.75) {};
	\end{tikzpicture}
	\begin{tikzpicture}[scale=1, every node/.style={scale=0.8}]
	\draw[dashed] (.5,.5) -- (3,3);
	\draw (.5,.5) -- (.5,1) -- (1,1) -- (1,1.5) -- (1.5,1.5) -- (1.5,2) -- (2,2) -- (2,2.5) -- (2.5,2.5) -- (2.5,3);
	\draw (.5,1.5) -- (1,1.5) -- (1,2) -- (1.5,2) -- (1.5,2.5) -- (2,2.5) -- (2,3);
	\draw[thick] (.5,.5) -- (.5,2) -- (1.5,2) -- (1.5,3) -- (3,3);
	\node[red] at (.75,.75) {$3$};
	\node[blue] at (.75,1.25) {$5$};
	\node[yellow] at (.75,1.75) {$6$};
	\node[blue] at (1.75,2.25) {$2$};
	\node[green] at (1.75,2.75) {$4$};
	\end{tikzpicture}
	\caption{A spanning tree with its corresponding sequence of shapes and parking function.}
	\label{fig:parkingoftree}
\end{figure}
	
It is well known that we can represent parking functions geometrically as Dyck paths with special labeling (see, e.g., \cite{Haglund2007}). To be precise, for a Dyck path $D$ of length $m$, we label the boxes to the right of north steps with integers $1$ through $m$ so that the numbers increase in columns upwards. The corresponding parking function $f_D:[m]\rightarrow [m]$ is then the map, where $f_D(i)$ is equal to the number of the column in which $i$ appears in $D$. As such, in \eqref{thm:formulawithtrees}, we exploit the notation and write $\mu(f)$ for the sequence of vertical-strip shapes corresponding to the Dyck path associated with $f$ (which is, in fact, a sequence of unicellular shapes since we have no diagonal steps in the path).
	
Also, recall that when we introduced Dyck paths in \cref{sec:schroderpaths}, we described a simple bijection between a certain subclass of Schr\"oder paths of length $m$ and Dyck paths of length $m-1$. In the complete graph case, we can combine it with \cref{prop:treestoschroder} to obtain a bijection $T\longleftrightarrow f$ satisfying $\nu(T) = \mu(f)$ with $T$ a Cayley tree on $m$ vertices and $f$ its corresponding parking function on $m-1$ cars (see \cref{fig:parkingoftree}). This assignment explains the second sum in \cref{thm:singcelldecomp}.
	
\begin{remark}
	Parking functions and Cayley trees are examples of classical combinatorial objects that appear in many different contexts. As such, finding bijections between the two classes that preserve certain statistics is an interesting problem in itself and the above correspondence is an example of that. Indeed, the bijection that naturally appears in this paper in the context of the combinatorics of symmetric functions, is the same as the one found recently in \cite[Section 7]{IrvingRattan2021}.
\end{remark}

\section*{Declarations}

\subsection*{Funding}

The author was supported by {\it Narodowe Centrum Nauki}, grant\linebreak UMO-2017/26/D/ST1/00186.

\subsection*{Conflicts of interests}

The author has no competing interests to declare that are relevant to the content of this article.

\bibliographystyle{amsalpha}

\bibliography{biblio2015}

\newcommand{\etalchar}[1]{$^{#1}$}
\def\cprime{$'$}
\providecommand{\bysame}{\leavevmode\hbox to3em{\hrulefill}\thinspace}
\providecommand{\MR}{\relax\ifhmode\unskip\space\fi MR }
\providecommand{\MRhref}[2]{%
  \href{http://www.ams.org/mathscinet-getitem?mr=#1}{#2}
}
\providecommand{\href}[2]{#2}
\begin{thebibliography}{HHL{\etalchar{+}}05b}

\bibitem[Ale21]{Alexandersson2021}
P.~Alexandersson, \emph{L{LT} polynomials, elementary symmetric functions and
  melting lollipops}, J. Algebraic Combin. \textbf{53} (2021), no.~2, 299--325.
  \MR{4238181}

\bibitem[AS22]{AlexanderssonSulzgruber2020}
Per Alexandersson and Robin Sulzgruber, \emph{A combinatorial expansion of
  vertical-strip llt polynomials in the basis of elementary symmetric
  functions}, Advances in Mathematics \textbf{400} (2022), 108256.

\bibitem[CM18]{CarlssonMellit2018}
E.~Carlsson and A.~Mellit, \emph{A proof of the shuffle conjecture}, J. Amer.
  Math. Soc. \textbf{31} (2018), no.~3, 661--697. \MR{3787405}

\bibitem[DF17]{DolegaFeray2017}
M.~Do{\l}{\k e}ga and V.~F\'{e}ray, \emph{Cumulants of {J}ack symmetric
  functions and the {$b$}-conjecture}, Trans. Amer. Math. Soc. \textbf{369}
  (2017), no.~12, 9015--9039. \MR{3710651}

\bibitem[DK21]{DolegaKowalski2021}
Maciej Dołęga and Maciej Kowalski, \emph{Llt cumulants and graph coloring}.

\bibitem[Do{\l}17]{Dolega2017}
M.~Do{\l}{\k{e}}ga, \emph{Strong factorization property of {M}acdonald
  polynomials and higher-order {M}acdonald's positivity conjecture}, J.
  Algebraic Combin. \textbf{46} (2017), no.~1, 135--163. \MR{3666415}

\bibitem[Do{\l}19]{Dolega2019}
M.~Do{\l}{\k e}ga, \emph{Macdonald cumulants, {$G$}-inversion polynomials and
  {$G$}-parking functions}, European J. Combin. \textbf{75} (2019), 172--194.
  \MR{3862962}

\bibitem[GH07]{GrojnowskiHaiman2007}
I.~Grojnowski and M.~Haiman, \emph{{Affine {H}ecke algebras and positivity of
  {LLT} and {M}acdonald polynomials}}, Preprint, 2007.

\bibitem[Hag07]{Haglund2007}
James Haglund, \emph{The q, t-catalan numbers and the space of diagonal
  harmonics : with an appendix on the combinatorics of macdonald polynomials},
  2007.

\bibitem[Hag16]{Haglund2016}
Jim Haglund, \emph{The combinatorics of knot invariants arising from the study
  of macdonald polynomials}, 04 2016.

\bibitem[Hai01]{Haiman2001}
M.~Haiman, \emph{Hilbert schemes, polygraphs and the {M}acdonald positivity
  conjecture}, J. Amer. Math. Soc. \textbf{14} (2001), no.~4, 941--1006.
  \MR{1839919}

\bibitem[HHL05a]{HaglundHaimanLoehr2005}
J.~Haglund, M.~Haiman, and N.~Loehr, \emph{A combinatorial formula for
  {M}acdonald polynomials}, J. Amer. Math. Soc. \textbf{18} (2005), no.~3,
  735--761. \MR{2138143}

\bibitem[HHL{\etalchar{+}}05b]{HaglundHaimanLoehrRemmelUlyanov2005}
J.~Haglund, M.~Haiman, N.~Loehr, J.~B. Remmel, and A.~Ulyanov, \emph{A
  combinatorial formula for the character of the diagonal coinvariants}, Duke
  Math. J. \textbf{126} (2005), no.~2, 195--232. \MR{2115257}

\bibitem[IR21]{IrvingRattan2021}
John Irving and Amarpreet Rattan, \emph{Trees, parking functions and
  factorizations of full cycles}, European Journal of Combinatorics \textbf{93}
  (2021), 103257.

\bibitem[Lee18]{Lee2018}
Seung~Jin Lee, \emph{Linear relations on llt polynomials and their k-schur
  positivity for k=2}, Journal of Algebraic Combinatorics \textbf{53} (2018),
  973 -- 990.

\bibitem[LLT97]{LascouxLeclercThibon1997}
A.~Lascoux, B.~Leclerc, and J.-Y. Thibon, \emph{Ribbon tableaux,
  {H}all-{L}ittlewood functions, quantum affine algebras, and unipotent
  varieties}, J. Math. Phys. \textbf{38} (1997), no.~2, 1041--1068.
  \MR{1434225}

\bibitem[LT00]{LeclercThibon2000}
B.~Leclerc and J.-Y. Thibon, \emph{Littlewood-{R}ichardson coefficients and
  {K}azhdan-{L}usztig polynomials}, Combinatorial methods in representation
  theory ({K}yoto, 1998), Adv. Stud. Pure Math., vol.~28, Kinokuniya, Tokyo,
  2000, pp.~155--220. \MR{1864481}

\bibitem[Mac88]{Macdonald1988}
I.~G. Macdonald, \emph{A new class of symmetric functions}, Publ. IRMA
  Strasbourg \textbf{372} (1988), 131--171.

\end{thebibliography}

\end{document}